\documentclass[10pt]{amsart}
\usepackage{amsfonts}
\usepackage{amscd}
\usepackage{epsf}
\usepackage{amssymb}
\usepackage{amsthm}
\usepackage{amsmath}
\usepackage{pb-diagram}
\usepackage{graphicx}
\def\noi{\noindent}

\def\rarrow{\rightarrow}

\def\xrarrow{\xrightarrow} %right arrow {label on top}

\def\noteq{\neq}

\def\<{\left<}
\def\>{\right>}

\def\onto{\twoheadrightarrow}

\def\ov{\overline}
\def\wt{\widetilde}

\DeclareMathOperator{\simp}{simp}
 \DeclareMathOperator{\Hom}{Hom}
 
\DeclareMathOperator{\sgn}{sgn}

\newcommand{\field}[1]{\mathbb{#1}}
\newcommand{\ZZ}{\ensuremath{{\field{Z}}}}

\newcommand{\QQ}{\ensuremath{{\field{Q}}}}

\def\l{\ell}
\def\ll{\lambda}

\newcommand{\A}{\ensuremath{{\mathcal{A}}}}

\newcommand{\F}{\ensuremath{{\mathcal{F}}}}
\newcommand{\G}{\ensuremath{{\mathcal{G}}}}

\newcommand{\N}{\ensuremath{{\mathcal{N}}}}

\newcommand{\Z}{\ensuremath{{\mathcal{Z}}}}
\def\thin{\F in}

\def\g{\gamma}
\def\d{\partial}

\def\e{\epsilon}
\def\f{\phi}
\def\k{\kappa}

\def\s{\sigma}
\def\Sig{\Sigma}

\def\Aut{{\rm Aut}}

\def\Fat{{\F at}}
\def\ktilde{\wt{\k}}
\def\ccc{\wt{c}_\G}

\def\Gam{\Gamma}

\newtheorem{thm}{Theorem}[section]
\newtheorem{lem}[thm]{Lemma}
\newtheorem{cor}[thm]{Corollary}
\newtheorem{prop}[thm]{Proposition}
\theoremstyle{definition}
\newtheorem{defn}[thm]{Definition}
\newtheorem{rem}[thm]{Remark}
\newtheorem{eg}[thm]{Example}
\newtheorem{conj}[thm]{Conjecture}

\title[Graph cohomology and Kontsevich cycles]
{Graph cohomology and Kontsevich cycles}

\begin{document}

\author{Kiyoshi Igusa}

\begin{abstract} We use the duality between compactly supported cohomology of
the associative graph complex and the cohomology of the mapping
class group to show that the duals of the Kontsevich cycles
$[W_\ll]$ correspond to polynomials in the Miller-Morita-Mumford
classes. We also compute the coefficients of the first two terms
of this polynomial. This extends the results of
\cite{[I:MMM_and_Witten]}, giving a more detailed answer to a
question of Kontsevich \cite{[Kontsevich:Airy]} and verifying more
of the conjectured formulas of Arbarello and Cornalba
\cite{[Arbarello-Cornalba:96]}.
\end{abstract}

%\address{Department of Mathematics, Brandeis University, Waltham, MA 02454}

%\email{igusa@brandeis.edu}


%57N05: Topology of E2, 2-manifolds
%55R40: Homology of classifying spaces, characteristic classes
%57M15: Relations with graph theory


\thanks{Partially supported by the National Science Foundation.}

%\date{\today}

\maketitle

%%%%%%%%%%%%%%%%%%%%%%%%%%%%%%%%%%%%%%%%%%%%%%%%%%%%%%%%%%%%%%%%%%%%%%%%%%%%%%%%%%%%%%%%

\section*{Introduction}

This paper explains the relationship between the Kontsevich cycles
in associative graph homology and the Miller-Morita-Mumford
classes in the cohomology of the mapping class group. We use a
version of the forested graph complex of Conant and Vogtmann
\cite{[CV02]} to go from the double dual of graph homology to the
cohomology of the mapping class group and we use our cyclic set
cocycle \cite{[I:MMM_and_Witten]} to evaluate the
Miller-Morita-Mumford classes on the Kontsevich cycles.

Graph homology was introduced by M. Kontsevich in
\cite{[Kontsevich:Formal-non-com]} and in
\cite{[Kontsevich:Feynman-diagrams]}. He constructed three graph
complexes which are called the ``Lie'', ``associative'' and
``commutative'' graph complexes. The case we are studying is the
associative case. This is the cohomology of ribbon graphs (graphs
with cyclic orderings of the half-edges incident to each vertex).
In \cite{[Kontsevich:Formal-non-com]},
\cite{[Kontsevich:Feynman-diagrams]} Kontsevich constructed
homology and cohomology classes in the associative case. He also
outlined a proof of the theorem that the rational finitely
supported cohomology of the associative graph complex is
isomorphic to the homology of the mapping class group. We use the
proof of this theorem given in \cite{[CV02]}.

\begin{thm}[Kontsevich]  $H_c^n(\G^\ast;\QQ)\cong
H_n(\coprod_{g,s}BM_{g}^s;\QQ)$ where $\G^\ast$ is the (reindexed)
associative graph complex.
\end{thm}

In this paper we use the category of ribbon graphs $\Fat$ which is
homotopy equivalent to the disjoint union of classifying spaces of
mapping class groups over all $s\geq1$ with $s\geq3$ when $g=0$.
(See \cite{[I:BookOne]} where this is shown to follow from
Culler-Vogtmann \cite{[Culler-Vogtmann-86]}.)
\[
    |\Fat|\simeq\coprod_{g,s}BM_{g}^s
\]
Therefore, the rational homology of $\Fat$ is isomorphic to the
rational finitely supported cohomology of the associative graph
complex. We give an explicit rational chain homotopy equivalence
between the cellular chain complex of this category and the graph
cohomology complex $\G_\ast$. Then we show that the dual
Kontsevich cycles in graph cohomology are (pull-backs of)
polynomials in the ``adjusted'' Miller-Morita-Mumford classes.
(They are adjusted by subtracting certain boundary classes.)

The main results of this paper were announced in
\cite{[I:MMM_and_Witten]} with short proofs. This paper gives more
detailed proofs, expresses them in the language of graph
cohomology and also extends these results to the next case. The
calculation at the end of the paper shows that for $n\noteq1$ we
have
\[
    [W_{n,1}^\ast]=3(-2)^{n+3}(2n+1)!!(\ktilde_n\ktilde_1-\ktilde_{n+1})-
    (-2)^{n+2}(2n+5)!!\ktilde_{n+1}.
\]
For $n=1$ we divide the right hand side by $2$.

In more detail the contents of this paper are as follows. In the
first section we review Kontsevich's definition of graph homology
using Conant and Vogtmann's formula for the Kontsevich orientation
of a graph. We define $\G_\ast$ to be the integral finitely
supported cohomology of the associative graph complex. Thus, e.g.,
$\G_0$ is the group of all integer valued functions $f$ on the set
of all isomorphism classes $[\Gam]$ of oriented trivalent ribbon
graphs so that $f[\Gam]=0$ for all but a finite number of $[\Gam]$
and $f[-\Gam]=-f[\Gam]$ where $-\Gam$ is $\Gam$ with the opposite
orientation. This complex has an augmentation map
\[
    \e:\G_0\to\QQ
\]
given by sending each dual generator $[\Gam]^\ast$ to
$\frac{o(\Gam)}{|\Aut(\Gam)|}$ where $o(\Gam)=\pm1$ depending on
the orientation of $\Gam$.

We define the \emph{integral subcomplex} $\G_\ast^\ZZ$ of
$\G_\ast$ to be the subcomplex generated by
\[
    \<\Gam\>=|\Aut(\Gam)|[\Gam].
\]
In $\G_0^\ZZ$ these elements have augmentation $\pm1$.

For any partition $\ll=(\ll_1,\ll_2,\cdots,\ll_r)$ of $n=\sum
\ll_i$ we define the \emph{dual Kontsevich cycle} $W_\ll^\ast$ to
be the homomorphism
\[
    W_\ll^\ast:\G_{2n}\to\QQ
\]
which sends each $[\Gam]^\ast$ in the Kontsevich cycle $W_\ll$ to
$\frac{o(\Gam)}{|\Aut(\Gam)|}$. Since
\[
    W_\ll^\ast\<\Gam\>=o(\Gam)=\pm1
\]
these are integral cocycles on $\G_\ast^\ZZ$. These cocycles were
considered by Kontsevich in \cite{[Kontsevich:Airy]}. They are the
Poincar\'{e} duals of certain strata of the moduli space of stable
curves.

In graph cohomology the dual Kontsevich cycles are linear
combinations of cocycles given by partition functions associated
with certain $1$-dimensional $A_\infty$ algebras. We give a
detailed account of this construction, essentially repeating what
Kontsevich says in \cite{[Kontsevich:Feynman-diagrams]} using the
Conant-Vogtmann definition of graph orientation.

For every ribbon graph $\Gam$ we next construct an acyclic
$\ZZ$-augmented chain complex $F_\ast(\Gam)$ over $\G_\ast$ so
that $F_\ast$ gives an acyclic carrier, i.e., a functor from the
category $\Fat$ of all ribbon graphs to the category of augmented
chain complexes over $\G_\ast$. We call it the \emph{forest
carrier}. This determines a chain map from the cellular chain
complex $C_\ast(\Fat)$ to $\G_\ast^\ZZ$ which is unique up to
homotopy. We show that this map is a rational homotopy equivalence
by constructing a rational inverse
\[
    \psi:\G_\ast^\ZZ\otimes\QQ=\G_\ast\otimes\QQ\to C_\ast(\Fat;\QQ).
\]
The chain map $\psi$ is defined by \emph{dual cells} $D(\Gam)$
modelled on the Poincar\'{e} duals of the Kontsevich cycles.

In the second section we discuss the Stasheff polyhedron. We use
the Conant-Vogtmann version of Kontsevich orientation to determine
the intrinsic orientation of the Stasheff polyhedron which
corresponds to the sign of the simplices in the dual cell
$D(\Gam)$. One of the main purposes of this is to justify the sign
convention used in \cite{[I:MMM_and_Witten]}. We also use this
discussion to prove that the forest carrier $F_\ast$ is acyclic as
claimed in the previous section.

In the third section we discuss the relationship between the
adjusted Miller-Morita-Mumford (MMM) classes $\ktilde_k$ and the
dual Kontsevich cycles. We review the formula for the adjusted MMM
classes given by the cyclic set cocycle and we show that the dual
Kontsevich cycles are polynomials in the adjusted MMM classes.
This is a detailed version of a one page argument in
\cite{[I:MMM_and_Witten]}.

The last section contains the calculation of the coefficients of
$[W_{n,1}^\ast]$ as a polynomial in the adjusted MMM classes. We
use the figures from section 2 which were drawn with this second
purpose in mind.

This paper started with a conversation with Karen Vogtmann about
graph homology. I should also thank Michael Kleber for some very
helpful discussions.

\begin{enumerate}
  \item Kontsevich cycles.
  \begin{enumerate}
    \item Category of ribbon graphs $\Fat$.
    \item Associative graph cohomology $\G_\ast$.
    \item Cocycles $W_{\ll}^\ast$ in graph cohomology.
    \item {Partition functions}.
    \item The forest carrier $F_\ast$.
    \item Dual cells.
\end{enumerate}
  \item Stasheff associahedra.
  \begin{enumerate}
    \item Stasheff polyhedron $K^n$.
    \item The category $\A_{n+3}$.
    \item Orientation of $K^n$.
    \item Orientation of $K^{odd}$.
    \item Proof of Proposition \ref{prop:forested graph complex is
acyclic}.
\end{enumerate}
  \item Miller-Morita-Mumford classes.
  \begin{enumerate}
    \item Cyclic set cocycle.
    \item Adjusted MMM classes in $H^{2k}(\G_\ast;\QQ)$.
    \item Cup products of adjusted MMM classes.
    \item Computing the numbers $b_{n_\ast}^{k_\ast}$.
    \item Kontsevich cycles in terms of MMM classes.
    \item Computing $a_\ll^\mu$.
\end{enumerate}
  \item Some computations.
  \begin{enumerate}
    \item {The degenerate case $n=0$}
    \item {Computation of $b_{n,1}^{n+1}$}
    \item {Conjectures}
\end{enumerate}
\end{enumerate}

\vfill\eject
%sec 0

\setcounter{section}{0}

\section{Kontsevich cycles}

\begin{enumerate}
    \item Category of ribbon graphs $\Fat$.
    \item Associative graph cohomology $\G_\ast$.
    \item Cocycles $W_{\ll}^\ast$ in graph cohomology.
    \item {Partition functions}.
    \item The forest carrier $F_\ast$.
    \item Dual cells.
\end{enumerate}

We review the basic definitions and give an explicit rational
homotopy equivalence between the finitely supported cohomology of
the associative graph complex and the cellular chain complex of
the category of ribbon graphs.

\subsection{Category of ribbon graphs $\Fat$.}

By a \emph{ribbon graph} (also known as \emph{fat graph}) we mean
a finite connected graph together with a cyclic ordering on the
half-edges incident to each vertex. We will use the following set
theoretic model for the objects in the category of graphs.

\begin{defn}\label{def of a graph}
Choose a fixed infinite set $\Omega$ which is disjoint from its
power set.\footnote{This occurs, e.g., if every element of
$\Omega$ is a set having greater cardinality than $\Omega$.} Then
by a \emph{graph} we mean a finite subset of $\Omega$ (the set of
\emph{half-edges}) together with two partitions of the set.
\begin{enumerate}
    \item A partition into pairs of half-edges which we call
    \emph{edges}.
    \item A partition into sets of cardinality (=\emph{valence}) $\geq3$ which
    we call \emph{vertices}.
\end{enumerate}
To avoid straying too far from conventional terminology we refer
to the elements of a vertex as \emph{incident half-edges}.
Equivalently, we define \emph{incident} to mean not disjoint.
\end{defn}

If $e=\{e^-,e^+\}$ is an edge in $\Gam$ then the vertices
$v_1,v_2$ incident to $e^-,e^+$ are the \emph{endpoints} of $e$.
If the endpoints are equal then $e$ is a \emph{loop}. If $e$ is
not a loop then we can collapse $e$ to a point forming a new graph
\[
    \Gam/e
\]
with one fewer edge, one fewer vertex and two fewer half-edges
than $\Gam$. Set theoretically, $\Gam/e$ is given by merging
$v_1,v_2$ and deleting $e^-,e^+$.

If $\Gam$ is a ribbon graph and $e$ is an edge in $\Gam$ which is
not a loop then $\Gam/e$ can be given the structure of a ribbon
graph in the obvious way by letting the new vertex be cyclically
ordered as
\[
    v_\ast=(h_1,\cdots,h_n,k_1,\cdots,k_m)
\]
if $v_1=(e^-,h_1,\cdots,h_n)$ and $v_2=(e^+,k_1,\cdots,k_m)$.

Morphisms of graphs and ribbon graphs can be given by collapsing
several edges to points and by isomorphisms. In other words,
certain subgraphs will be collapsed to points.

By a \emph{subgraph} of a graph $\Gam$ we mean a subset of the set
of vertices together with all incident half edges and a set of
edges both endpoints of which lie in the chosen set of vertices.
For example, we could take all of the vertices and none of the
edges. A subgraph will usually not be a graph since it usually has
unpaired half-edges. The unpaired half-edges of a subgraph will be
called its \emph{leaves}. If the graph $\Gam$ is
connected\footnote{A graph/subgraph is \emph{connected} if it is
not the disjoint union of two graphs/subgraphs.} then every
subgraph is determined by its set of leaves.

A subgraph is a \emph{tree} if it is connected and has one more
vertex than edge. A \emph{forest} is a disjoint union of trees. If
$F$ is a spanning\footnote{\emph{Spanning} means containing all
the vertices.} forest in a graph $\Gam$, let $\Gam/F$ be the graph
obtained by collapsing each tree in $F$ to a separate point. By
this we mean collapse the edges of each tree to a point. Thus:
\begin{enumerate}
    \item The edges of $\Gam/F$ are the edges of $\Gam$ which do
    not lie in $F$.
    \item The vertices of $\Gam/F$ are the sets of leaves of the
    component trees of $F$.
\end{enumerate}

\begin{defn}\label{def of graph morphisms}
A \emph{morphism} of graphs $\f:\Gam_0\to\Gam_1$ is defined to be
an isomorphism \begin{equation}\label{morphism of graphs}
    \Gam_0/F\cong\Gam_1
\end{equation}
for some spanning forest $F$ in $\Gam_0$. In other words, the
inverse image of every edge in $\Gam_1$ is an edge in $\Gam_0$ and
the inverse image of every vertex of $\Gam_1$ is a tree in
$\Gam_0$.
\end{defn}

One thing is obvious from this definition. A morphism
$\f:\Gam_0\to\Gam_1$ is uniquely determined by the value of
$\f^{-1}(e)$ for every half edge $e$ in $\Gam_1$. The reason is
that this information specifies the forest $F$ and also gives an
isomorphism $\Gam_1\cong\Gam_0/F $.

Morphisms of graphs (and ribbon graphs) also have the following
left cancellation property.

\begin{prop}\label{prop:left cancellation}
Any two morphisms $f,g:\Gam_0\to \Gam_1$ which are equalized by a
morphism $h:\Gam_1\to\Gam_2$ are equal. I.e., \[hf=hg\
\Rightarrow\ f=g.\]
\end{prop}

\begin{rem}\label{rem:morphisms are mono and epi}
In category theoretic terminology, this proposition says that
morphisms of graphs are \emph{monomorphisms}. They are also
obviously \emph{epimorphisms}. I.e., they satisfy both left and
right cancellation.
\end{rem}

\begin{proof}
In order for $f,g$ to be different, there must be a half-edge $e$
in $\Gam_1$ so that $f^{-1}(e)\noteq g^{-1}(e)$. But $hf=hg$
cannot send two different half-edges of $\Gam_0$ to the same
half-edge in $\Gam_2$. So $h(e)$ must be a vertex $v$. The inverse
image of $v$ is a tree $T_0$ in $\Gam_0$ and another tree $T_1$ in
$\Gam_1$ and the leaves of both trees map bijectively onto the
half-edges incident to $v$. Consequently, $f,g$ give the same
bijection of the leaves of $T_0$ with the leaves of $T_1$. Any
edge in $T_0$ is uniquely characterized by the partitioning of the
set of leaves which would result if we cut the edge. And each
interior half-edge of $T_0$ is determined by the corresponding
subset of the set of leaves. Thus $f^{-1}(e)=g^{-1}(e)$ which is a
contradiction.
\end{proof}

If $\Gam$ is a ribbon graph then the set of leaves of every tree
in $\Gam$ inherits a cyclic order. Consequently, $\Gam/F$ has an
induced structure as a ribbon graph. A graph morphism
$\f:\Gam_0\to\Gam_1$ will be called a \emph{ribbon graph morphism}
if it respects these cyclic orderings, i.e., if (\ref{morphism of
graphs}) is an isomorphism of ribbon graphs.

If $\Gam,\Gam'$ are ribbon graphs, let $\Hom(\Gam,\Gam')$ denote
the set of all ribbon graphs morphisms $\Gam\to\Gam'$. Since this
is a subset of the set of all graph morphisms, left and right
cancellation hold for these morphisms as well. Thus we get the
following corollary where $\Aut(\Gam)=\Hom(\Gam,\Gam)$ is the
group of ribbon graph automorphisms of $\Gam$.

\begin{cor}\label{cor:Aut(G) acts freely on Hom(G,G')}
$\Aut(\Gam)$ acts freely on the right on $\Hom(\Gam,\Gam')$ and
$\Aut(\Gam')$ acts freely on the left.
\end{cor}

Let $\Fat$ denote the category of all ribbon graphs and ribbon
graph morphisms. Let $|\Fat|$ denote its geometric realization:
\[
    |\Fat|=\coprod_{n}\coprod_{\Gam_\ast\in\N_n\Fat}\Delta^n/\sim
\]
Then we have the following theorem which I learned from Penner
\cite{[Penner87:decorated]}, which goes back to Strebel
\cite{[Strebel]} but which I prove using Culler-Vogtmann
\cite{[Culler-Vogtmann-86]}. For details, see \cite{[I:BookOne]}.

\begin{thm}
\[
    |\Fat|\simeq\coprod_{g,s} BM_g^s
\]
where $M_g^s$ is the \emph{mapping class group} of a surface of
genus $g$ with $s$ punctures.\footnote{$M_g^s$ is the group of
homotopy classes of orientation preserving self homeomorphisms of
a connected Riemann surface of genus $g$ with $s$ punctures. This
group maps onto the symmetric group on $s$ letters with kernel
equal to the mapping class group $M_{g,s}$ of genus $g$ surfaces
with $s$ marked points.}
\end{thm}

Many ribbon graphs have an intrinsic orientation in the sense of
Kontsevich graph homology.

\subsection{Associative graph cohomology $\G_\ast$.}

Graph homology (of ribbon graphs) is rationally dual to the
homology of the category of ribbon graphs. More precisely, we have
an isomorphism between rational compactly supported cohomology of
the associative graph complex $\G^\ast$ and the rational homology
of the mapping class group. We will construct an integral chain
map which realizes this rational equivalence.

However, the main purpose of introducing graph homology in the
present context is to fix our orientation conventions. We use the
definitions given in \cite{[CV02]}. Since our graphs are connected
this agrees with Kontsevich's orientation convention.

\begin{defn}[Conant-Vogtmann]\label{def:orientation of a graph}
An \emph{orientation} of a graph is defined to be an orientation
on the vector space spanned by the set of vertices and half-edges.
\end{defn}

This means that an orientation of $\Gam$ can be given by ordering
the vertices and orienting the edges of $\Gam$.

\begin{rem}[Conant-Vogtmann]\label{def:natural orientation of odd valent
ribbon graphs} A ribbon graph has a natural orientation if all its
vertices have odd valence. Since such a graph necessarily has an
even number of vertices it does not matter if we put the vertex
first and then the incident half-edges in cyclic order or the
other way around.
\end{rem}

\begin{defn}[Conant-Vogtmann]\label{def:induced orientation of Gam
mod e} Given a graph $\Gam$ and an edge $e$ in $\Gam$ which is not
a loop, let $\Gam/e$ be the graph obtained from $\Gam$ by
collapsing the edge $e$. If an orientation on $\Gam$ is given by
orienting all edges and ordering the vertices so that the source
of $e$ is first and its target is second, the \emph{induced
orientation} on $\Gam/e$ is given by taking the coalesced vertex
to be first and letting the remaining vertices and edges be
ordered and oriented as before.
\end{defn}

\begin{rem}\label{rem:rule for collapsing edges}
If $e$ is an oriented edge going from vertex $v_i$ to vertex $v_j$
then, when we collapse $e$, we should insert the new vertex in
position $i$ and multiply the orientation by $(-1)^j$. This will
give the induced orientation for $\Gam/e$.
\end{rem}

The \emph{associative graph homology complex} can now be defined.
For all $n\geq0$ let $\G^n$ be the free abelian group generated by
all isomorphism classes $[\Gam]$ of (connected) oriented ribbon
graphs $\Gam$ of codimension\footnote{A graph has
\emph{codimension} $n$ if it is obtained from a trivalent graph by
collapsing $n$ edges.} $n$ module the relation $-[\Gam]=[-\Gam]$
where $-\Gam$ is $\Gam$ with the opposite orientation. If $\Gam$
has an orientation reversing automorphism this implies that
$2[\Gam]=0$. Define the boundary operator $\d:\G^n\to\G^{n+1}$ by
\[
    \d[\Gam]=\sum_e [\Gam/e]
\]
where the sum is over all edges in $\Gam$ which are not loops.

The compactly supported dual of this complex is the
\emph{(associative) graph cohomology complex} given as follows.

\begin{defn}
For all $n\geq0$ let $\G_n$ be the additive group of all
homomorphisms $f:\G^n\to\ZZ$ so that $f[\Gam]\noteq0$ for only
finitely many $[\Gam]$. (In particular, $f[\Gam]=0$ if $\Gam$ has
an orientation reversing automorphism.) Thus $\G_n$ is generated
by duals $[\Gam]^\ast$ of generators of $\G^n$. The boundary map
$d:\G_n\to \G_{n-1}$ is given in terms of these dual generators by
\begin{equation}\label{eq:boundary in Gn}
    d[\Gam]^\ast=\sum \l_i[\Gam_i]^\ast
\end{equation}
where $\l_i$ is equal to the number of edges $e$ in $\Gam_i$ so
that $\Gam_i/e\cong\Gam$ minus the number of edges in $\Gam_i$ so
that $\Gam_i/e\cong-\Gam$. The sum is over a basis for $\G_{n-1}$.
\end{defn}

The coefficient $\l_i$ in (\ref{eq:boundary in Gn}) can be written
as
\[
    \l_i=\frac{|\Hom^+(\Gam_i,\Gam)|-|\Hom^-(\Gam_i,\Gam)|}
    {|\Aut(\Gam)|}\in\ZZ
\]
where $\Hom^\pm(\Gam_i,\Gam)$ is the set of morphisms
$f:\Gam_i\to\Gam$ so that the orientation of $\Gam$
agrees/disagrees with the orientation induced from $\Gam_i$ by
$f$. In other words, $\l_i$ is the number of left equivalence
classes of morphisms $\Gam_i\to\Gam$ counted with sign.

Let $r_i$ be the number of right equivalence classes of such maps
counted with sign. Then
\[
    r_i=\frac{|\Hom^+(\Gam_i,\Gam)|-|\Hom^-(\Gam_i,\Gam)|}
    {|\Aut(\Gam_i)|}\in\ZZ.
\]
So (\ref{eq:boundary in Gn}) can be written as
\[
    d\<\Gam\>=\sum r_i\<\Gam_i\>
\]
where
\[
    \<\Gam\>:=|\Aut(\Gam)|[\Gam]^\ast.
\]

\begin{defn}
Let $\G_\ast^\ZZ$ denote the subcomplex of $\G_\ast$ generated by
the elements $\<\Gam\>$. We call $\G_\ast^\ZZ$ the \emph{integral
subcomplex} of $\G_\ast$.
\end{defn}

The boundary map in $\G_\ast^\ZZ$ can be described in terms of
expanding vertices. If $\Gam$ is an oriented ribbon graph, each
vertex of valence $n$ can be expanded into two vertices connected
by an edge in
\[
    \frac{n^2-3n}{2}
\]
different ways. Each of these choices gives a ribbon graph
$\Gam_i$ with a distinguished edge $e$ and an isomorphism
\[
    \Gam_i/e\cong\Gam
\]
so that $\Gam_i$ is unique up to isomorphism over $\Gam$. We give
$\Gam_i$ the orientation induced from $\Gam$ by this isomorphism.
The boundary map of the integral subcomplex is then given by
\[
    d\<\Gam\>=\sum \<\Gam_i\>.
\]

\begin{eg}\label{ex:boundary map in AG}
Consider the ribbon graphs $\Gam_0,\Gam_0',\Gam_1$ shown in Figure
\ref{fig:D8 graph}. We take the natural orientation on the
trivalent graphs $\Gam_0, \Gam_0'$ and the induced orientation on
$\Gam_1$ given by the isomorphism
\[
    \Gam_1\cong\Gam_0/e.
\]
If we write the orientation of $\Gam_0$ as
\[
    v_1e_1bc\,v_2e_2da\cdots
\]
then we see that the induced orientation for $\Gam_1$ is
$vabcd\cdots$. The other three straight edges of $\Gam_0$ are
equivalent to $e$ since $\Aut(\Gam_0)\cong D_8$ acts transitively
on the set of $4$ straight edges. So the boundary of $[\Gam_0]$ in
the associative graph complex is given by
\[
    \d[\Gam_0]=4[\Gam_1]+8[\Gam_1']
\]
where $\Gam_1'$ is given by collapsing one of the curved edges of
$\Gam_0$.
%*
\begin{figure}
\includegraphics{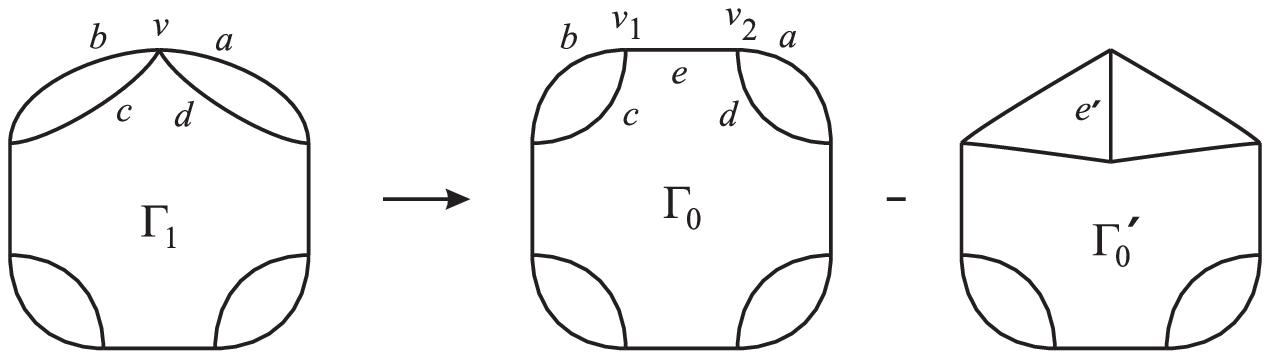}
\caption{$d\<\Gam_1\>=\<\Gam_0\>-\<\Gam_0'\>$.}\label{fig:D8
graph}
\end{figure}
%*

The ribbon graph $\Gam_1$ is also part of the boundary of
$[\Gam_0']$ since
\[
    -\Gam_1\cong\Gam_0'/e'.
\]
But there are no other edges in $\Gam_0'$ equivalent to $e'$.
Consequently,
\[
    d[\Gam_1]^\ast=4[\Gam_0]^\ast-[\Gam_0'].
\]
Since the orders of the automorphism groups are: $2,8,2$,
respectively, we get:
\[
    d\<\Gam_1\>=\<\Gam_0\>-\<\Gam_0'\>.
\]
\end{eg}

We will be looking at the rational cochain complex
\[
    \Hom(\G_\ast,\QQ).
\]
This is the rational double dual of the original graph homology
complex $\G^\ast$. Thus cocycles in this complex, such as
$W_{\ll}^\ast$ defined below, are ``infinite cycles'' in the graph
homology complex.

\begin{rem} Since $\G_n^\ZZ$ is a free abelian group whose
generators $\<\Gam\>$ form a $\QQ$-basis for $\G_n\otimes\QQ$, its
integral dual forms a lattice
\[
    \Hom(\G_\ast^\ZZ,\ZZ)\subseteq\Hom(\G_\ast,\QQ)
\]
which we call the \emph{integral cochain complex}. Elements of
this subcomplex will be called \emph{integral cochains} on
$\G_\ast$.
\end{rem}

\subsection{Cocycles $W_{\ll}^\ast$ in graph cohomology.}

\begin{defn}
If $\ll=(\ll_1,\ll_2,\cdots,\ll_r)$ is a sequence of positive
integers let $W_{\ll}$ be the set of all ribbon graphs $\Gam$
which are trivalent at all but $r$ vertices $v_1,\cdots,v_r$ which
have valence $2\ll_i+3$, resp. This set will be called the
\emph{Kontsevich cycle}. The \emph{dual Kontsevich cycle}
$W_{\ll}^\ast\in \G^{2|\ll|}$ (where $|\ll|=\ll_1+\ll_2+\cdots
+\ll_r$) is given by
\[
    W_{\ll}^\ast[\Gam]^\ast:=
    \begin{cases}
    \frac{o(\Gam)}{|\Aut(\Gam)|} &\text{if  }\Gam\in W_{\ll}\\
    0 &\text{if  }\Gam\notin W_{\ll}
    \end{cases}
\]
where $o(\Gam)=\pm1$ depending on whether the given orientation of
$\Gam$ agrees with the natural orientation (Remark
\ref{def:natural orientation of odd valent ribbon
graphs}).\end{defn}

The dual Kontsevich cycle are integral cochains since they can be
given by
\[
    W_\ll^\ast\<\Gam\>=o(\Gam)
\]
if $\Gam\in W_\ll$. Also, note that $W_\ll, W_\ll^\ast$ are
independent of the order of the $\ll_i$.

In the special case when $r=1$, $W_{k}$ is called the \emph{Witten
cycle} and $W_{k}^\ast$ will be called the \emph{dual Witten
cycle}. In the case $r=0$, $W_\emptyset$ is the set of all
trivalent (connected) ribbon graphs and
\[
    W_\emptyset^\ast=\e:\G_0\to\QQ
\]
is the map sending the dual $[\Gam]^\ast$ of every trivalent graph
$\Gam$ to $o(\Gam)/|\Aut(\Gam)|$. On the integral subcomplex this
gives an epimorphism
\[
    \e:\G_0^\ZZ\onto\ZZ
\]
sending each generator $\<\Gam\>$ to $o(\Gam)=\pm1$. We define
these maps to be the \emph{augmentation maps} for $\G_\ast$ and
$\G_\ast^\ZZ$.

We will also consider \emph{degenerate} cases where some of the
indices are zero. We interpret these $0$'s as counting the number
of distinct trivalent vertices:
\[
    W_{0^k}^\ast[\Gam]^\ast:=
        \begin{cases}
    \binom{n}{k}\frac{o(\Gam)}{|\Aut(\Gam)|} &\text{if  }\Gam\in W_{\emptyset}\text{ with }n\text{ vertices}\\
    0 &\text{if  }\Gam\text{ is not trivalent.}
    \end{cases}
\]

\begin{prop}\label{prop:Xk is a cocycle}
Each dual Kontsevich cycle $W_{\ll}^\ast$ is an integral cocycle.
\end{prop}

\begin{proof}
Let $n=2|\ll|$. Then we want to show that
\[
    W_{\ll}^\ast(d\<\Gam\>)=0
\]
for all oriented ribbon graphs $\Gam$ of codimension $n+1$.
However, the only case in question occurs when $\Gam$ has only one
even valent vertex, call it $v_0$. The orientation on $\Gam$ can
be given by first taking $v_0$, then the incident half-edges
$e_1,\cdots,e_{2m}$, then all other vertices with their incident
half-edges in cyclic order. The orientation depends on which of
the half-edges at $v_0$ is first.

There are three cases.

\underline{Case 1}. The graph $\Gam$ has $r-2$ odd valent vertices
of codimension\footnote{The \emph{codimension} of a vertex is
equal to its valence minus $3$.} $\geq2$. After re-indexing the
$\ll_i$ we may assume that these codimensions are
$2\ll_3,2\ll_4,\cdots,2\ll_r$. The even valent vertex $v_0$ must
have valence $2\ll_1+2\ll_2+4$ and it needs to split into two
vertices of codimension $2\ll_1$ and $2\ll_2$. There is always an
even number of ways to do this ($2\ll_1+2$ ways if $\ll_1=\ll_2$
and $2\ll_1+2\ll_2+4$ ways if not) and half of them will give one
sign and half the other. (The sign alternates as we rotate the
half edges incident to $v_0$.) Consequently, the value of
$W_{\ll}^\ast$ on $d[\Gam]^\ast$ will be zero.

\underline{Case 2}. The graph $\Gam$ has $r-1$ odd valent vertices
of codimension $\geq2$. We may assume that these codimensions are
$2\ll_2,\cdots,2\ll_r$ (after re-indexing the $\ll_i$). The vertex
$v_0$ must have valence $2\ll_1+4$ and it needs to split into two
vertices of valence $2\ll_1+3$ and $3$. There are $2\ll_1+4$ ways
to do this and half of them will give one sign and half the other.

\underline{Case 3}. $v_0$ has valence $4$. It can split into two
trivalent vertices in two ways with opposite sign as we saw in
Example \ref{ex:boundary map in AG}.
\end{proof}

Proposition \ref{prop:Xk is a cocycle} also follows from an
observation of Kontsevich that $A_\infty$ superalgebras give
partition functions on ribbon graphs which are cocycles on
associative graph cohomology.

\subsection{Partition functions}

In \cite{[Kontsevich:Formal-non-com]} and
\cite{[Kontsevich:Feynman-diagrams]}, Kontsevich explains how a
finite dimensional $A_\infty$ superalgebra $A$ gives a cocycle on
the associative graph cohomology complex $\G_\ast$. Kontsevich
assumed that $A$ was an algebra over the real numbers. However, it
is easy to see that the ground field can have any characteristic.
In fact, we only need to assume that $A$ is a finitely generated
free module over a commutative ring $R$.

We will go over the definition of an $A_\infty$ superalgebra
following Getzler and Jones \cite{[GetzlerJones89]}. Then we
revise Kontsevich's definition of the partition function using the
Conant-Vogtmann definition of graph orientation. Finally, we
examine the special case of one dimensional algebras to verify
Kontsevich's claim \cite{[Kontsevich:Feynman-diagrams]} that the
cocycles coming from these examples linearly span the space of
polynomials in the Miller-Morita-Mumford classes. Translated into
the present setting these cocycles are easily seen to be linear
combinations of the dual Kontsevich cycles (which come from
Kontsevich's earlier paper \cite{[Kontsevich:Airy]}).

\begin{defn}\label{def:A infinity superalgebra}
By an \emph{$A_\infty$ superalgebra} we mean a $\ZZ/2$-graded
algebra $A=A_0\oplus A_1$ over a commutative ring $R$ together
with a sequence of $R$-linear mappings
\[
    m_k:A^{\otimes k}\to A,\quad k\geq 1
\]
which are homogeneous of degree $k$ (mod $2$) so that for
homogeneous elements $x_1,\cdots,x_k$ we have
\[
    \sum_{r+s+t=k}(-1)^um_{r+1+t}(x_1,\cdots,x_r,
    m_s(x_{r+1},\cdots,x_{r+s}),x_{r+s+1},\cdots,x_k)=0
\]
where $u={r+st+s|x_1|+\cdots+s|x_r|}$.
\end{defn}

Suppose that $A\cong R^n$ is finitely generated and free as an
$R$-module. Suppose that $m_1=0$. And suppose that we have an
\emph{nondegenerate even scalar product}
\[
    \<\ , \>:A\otimes A\to R.
\]
This means the following.
\begin{enumerate}
    \item $\<a,b\>=0$ if $|a|+|b|=1$.
    \item $\<a,b\>=(-1)^{|a|}\<b,a\>$. (This implies (1) if
    $2$ is not a zero divisor in $R$.)
    \item There is a degree $0$ $R$-linear isomorphism
    $$D:\Hom_R(A,R)\xrarrow{\approx}
    A$$ so that $\<a,Df\>=f(a)$.
    \item For all $x_0,\cdots,x_n\in A_0\coprod A_1$ we
    have
    \[
\<m_n(x_1,\cdots,x_n),x_0\>=
(-1)^{n+|x_0|+n|x_0|}\<m_n(x_0,\cdots,x_{n-1}),x_n\>
    \]
\end{enumerate}
Then we have a \emph{partition function}
\[
    Z_A:\G_\ast^\ZZ\to R
\]
given on any generator $\<\Gam\>$ as follows.

First, choose an ordering for the vertices $v_1,v_2,\cdots$ of
$\Gam$. Next, label the half-edges incident to each $v_i$ in
\emph{reverse (clockwise)} order $e_{i1},\cdots,e_{in_i},e_{i0}$
(if $v_i$ has valence $n_i+1$). Let $\e_1=\pm1$ so that
\[
    o(\Gam)=\e_1\sgn(v_1,e_{10},e_{1n_1},\cdots,
    e_{11},v_2,e_{20},e_{2n_2},\cdots,
    e_{21},v_3,\cdots).
\]
Choose an $R$-basis $b_1,\cdots,b_n$ for $A$ ($b_i\in A_0\coprod
A_1$) and a dual basis $b_1^\ast,\cdots,b_n^\ast\in \Hom_R(A,R)$
so that
\[
    \<b_i,Db_j^\ast\>=b_j^\ast(b_i)=\delta_{ij}.
\]
Then the partition function is given by the state sum
\begin{equation}\label{eq:partition function of A}
    Z_A\<\Gam\>=\e_1\sum_{states}
    \prod_i\<m_{n_i}(x_{i1},\cdots,x_{in_i}),x_{i0}\>
    \e_2\prod_j\<D\ov{y}_j^\ast,D{y}_j^\ast\>
\end{equation}
The sum is over all {states} where a \emph{state} of $\Gam$ is
given by assigning a basis element $x_{ij}$ of $A$ to each half
edge $e_{ij}$ of $\Gam$. The first product is over all vertices
$v_i$. The second product is over all edges $(h_j,\ov{h}_j)$. Here
$y_j^\ast$ represents the dual basis element corresponding to the
basis element $y_j$ assigned to $h_j$ and similarly for
$\ov{y}_j^\ast$. The sign $\e_2=\pm1$ is the sign of the
permutation of the odd half-edges (those assigned elements of
$A_1$ as basis elements) as they appear in the sequence:
\begin{equation*}
    e_{11},e_{12},\cdots,e_{1n_1},e_{10},e_{21},\cdots,e_{2n_2},e_{20},e_{31},\cdots
\end{equation*}
which places each next to its other half (placing $\ov{h}_j$ next
to and on the right of $h_j$).

\begin{thm}[Kontsevich]\label{thm:ZA is a cocycle}
$Z_A$ is a cocycle on $\G_\ast^\ZZ$.
\end{thm}

\begin{rem}
The usual definition of the partition function has a factor of
$\frac1{|\Aut(\Gam)|}$:
\[
    Z_A[\Gam]^\ast=\frac{\e_1}{|\Aut(\Gam)|}\sum_{states}
    \prod_i\<m_{n_i}(x_{i1},\cdots,x_{in_i}),x_{i0}\>
    \e_2\prod_j\<D\ov{y}_j^\ast,D{y}_j^\ast\>
\]
This factor disappears on the integral subcomplex $\G_\ast^\ZZ$
since $\<\Gam\>={|\Aut(\Gam)|}[\Gam]^\ast$.
\end{rem}

\begin{proof} It is easy to see that the partition function is
well-defined. Since each basis element $a_{ij}$ and its dual
appears exactly once, the formula (\ref{eq:partition function of
A}) for $Z_A(\Gam)$ is independent of the choice of basis. If we
transpose two vertices $v_1,v_2$ then both $\e_1$ and $\e_2$
change by a factor of $(-1)^{n_1n_2}$. Finally, if we cyclically
permute the half-edges around a vertex $v$ of valence $n+1$ then
the signs $\e_1,\e_2$ and the value change by factors of
\begin{align*}
    \e_1'/\e_1&=(-1)^n\\
    \e_2'/\e_2&=(-1)^{|x_0|\sum|x_i|}\\
    value'/ value&=(-1)^{n+|x_0|+|x_0|n}
\end{align*}
The product of these factors is $1$ since the degrees of
$x_0,\cdots,x_n$ must add up to $n$ mod $2$. (Otherwise the
expression (\ref{eq:partition function of A}) is zero.) Thus $Z_A$
is well-defined. It remains to show that $Z_A$ is a cocycle.

The boundary of any generator $\<\Gam\>$ in $\G_\ast^\ZZ$ is a sum
over all vertices $v$ of $\Gam$ of all ribbon graphs $\Gam'$
obtained from $\Gam$ by expanding $v$ into two vertices. For each
fixed $v$ the sum of the values of $Z_A(\Gam')$ add up to zero. To
see this we label the half-edges clockwise around $v$. This means
that the Conant-Vogtmann orientation starts as:
\[
    o(\Gam)=\sgn(v\, e_0\, e_n \cdots  e_1\, v' \cdots).
\]
When we expand $v$ we get $\Gam'$ with orientation
\[
    o(\Gam')=\sgn( v_1\,v_2\,h\,\ov{h}\,e_0\,e_n\cdots e_1\,v'\cdots)
\]
\[
    =(-1)^u\sgn(v_1\,h\,e_{r+s}\cdots e_{r+1}\,
    v_2\,e_0\,e_n\cdots e_{r+s+1}\,\ov{h}\,e_{r+1}\cdots e_1\,v'\cdots)
\]
where $n=r+s+t$ with
\[
    \e_1=(-1)^u=(-1)^{st+s+t+1}=(-1)^{r+st+n+1}
\]
The corresponding terms of the partition function are:
\begin{equation}\label{eq:used in proof that ZA is a cocycle}
    \<m_s(x_{r+1},\dots, x_{r+s}),y\>
    \<m_{r+t+1}(x_1,\dots, x_r,\ov{y},x_{r+s+1},\dots,x_n),x_0\>
    \<D\ov{y}^\ast,Dy^\ast\>
\end{equation}
with associated relative sign term
\[
    \e_2=(-1)^{s|x_1|+\cdots+s|x_r|}
\]
since the degrees of $x_{r+1},\cdots,x_{r+s},y$ must add up to
$r$.

Using the identity
\[
    \sum_i\<x,b_i\>\<y,Db_i^\ast\>=\<x,y\>
\]
we see that the expression (\ref{eq:used in proof that ZA is a
cocycle}) contracts to
\[
    \<m_{r+t+1}(x_1,\dots, x_r,m_s(x_{r+1},\dots, x_{r+s}),x_{r+s+1},\dots,x_n),x_0\>
\]
when summed over all allowed values of
$y,\ov{y},y^\ast,\ov{y}^\ast$. By definition of an $A_\infty$
algebra, the product of this with $\e_1\e_2$ adds up to zero if we
sum over all $\Gam'$ obtained from $\Gam$ by expanding $v$ since
$n$ is constant. We need the assumption $m_1=0$ since $\Gam'$ has
no bivalent vertices.
\end{proof}

Suppose that $x=(x_0,x_1,x_2,\cdots)$ is an infinite sequence of
rational numbers. {Then Kontsevich points out that there is a
$1$-dimensional $A_\infty$ algebra $A=A_0=\QQ$ with scalar product
in which $m_{2k}$ is multiplication by $x_{k-1}$, $m_{odd}=0$ and
$\<a,b\>=ab$. Since the states of $\Gam$ are given by assigning a
basis vector to each half-edge, there is only one state and the
partition function $Z_x=Z_A$, which is a sum over all states, has
only one term.}

\begin{eg}\label{eg:partition function Zx}
The \emph{partition function}
\[
    Z_x:\G_\ast\to \QQ
\]
is the cocycle defined by the equation
\[
    Z_x[\Gam]^\ast=\frac{o(\Gam)}{|\Aut(\Gam)|}{x_0^{r_0}x_1^{r_1}\cdots}
\]
if $\Gam$ is a ribbon graph with $r_i$ vertices of valence $2i+3$
for $i=0,1,2,\cdots$ and no vertices of even valence.
\end{eg}

Since the Euler characteristic of $\Gam$ is given by
\[
    \chi(\Gam)=-\frac12\sum r_i(2i+1),
\]
the value of $r_0$ can be written as
\begin{equation}\label{eq:value of r0}
    r_0=-2\chi-\sum_{i\geq1}r_i(2i+1).
\end{equation}
Thus, the partition function $Z_x$ can be given in terms of the
dual Kontsevich cycles by
\begin{equation}\label{eq:formula for Zx in terms of W ll}
    Z_x=x_0^{-2\chi}\sum_\ll y^\ll W_\ll^\ast
\end{equation}
where $y^\ll=\prod_i(x_i/x_0^{2i+1})^{r_i}$ if
$\ll=(1^{r_1},2^{r_2},\cdots)$ and the sum is over all $\ll$ so
that $r_0$, as given by (\ref{eq:value of r0}), is nonnegative.
(So the right hand side of (\ref{eq:formula for Zx in terms of W
ll}) is well-defined only when $\chi$ is fixed.)

Thus, if we restrict to the subcomplex of the graph cohomology
complex $\G_\ast$ generated by $[\Gam]^\ast$ where $\chi(\Gam)$ is
fixed, the linear span of these partition functions is the same as
the linear span of the $W_\ll^\ast$ and, by Corollary
\ref{cor:span W ll = QQ[k1,k2,etc]} below, this is equal to the
algebra generated by the adjusted Miller-Morita-Mumford classes
$\ktilde_k$ excluding $\ktilde_0=\chi$. (However, we lose the
linear independence of the $[W_\ll^\ast]$ when we restrict to this
finitely generated subcomplex of $\G_\ast$.)

Setting $x_0=1$ and taking partial derivatives of (\ref{eq:formula
for Zx in terms of W ll}) with respect to $y^\ll$ using finite
differences, we can conclude that, for each fixed $\chi$, the dual
Kontsevich cycles $W_\ll^\ast$ are linear combinations of
partition functions $Z_x$ for various multi-indices $x$.
Consequently, Kontsevich's theorem (\ref{thm:ZA is a cocycle})
that each $Z_x$ is a cocycle implies that each $W_\ll^\ast$ is a
cocycle.

To pull the cocycles $W_{\ll}^\ast$ back to the category of ribbon
graphs we need to use an acyclic carrier related to the forested
graph complex of \cite{[CV02]}.

\subsection{The forest carrier $F_\ast$.}

Conant and Vogtmann use forested graph complexes to show that
graph homology is rationally isomorphic to the cohomology of the
mapping class group as claimed by Kontsevich. In our notation the
forested graph complex is the total complex of an integral acyclic
carrier
\[
    F_\ast:\Fat\to \G_\ast^\ZZ\subseteq\G_\ast
\]
which we call the ``forest carrier.''

Suppose that $\Gam_0$ is a ribbon graph. Then we will construct a
chain complex $F_\ast(\Gam_0)$ generated by the isomorphism
classes of all ribbon graphs $\Gam$ which map to $\Gam_0$. Each
such object is given by a ``{forested ribbon graph},'' i.e., a
ribbon graph $\Gam$ with a spanning forest $F$ (so that
$\Gam/F\cong\Gam_0$). If $\Gam_0$ has codimension $n$ then
$F_\ast(\Gam_0)$ will be the augmented chain complex
\[
    0\to F_n(\Gam_0)\to F_{n-1}(\Gam_0)\to\cdots\to
    F_0(\Gam_0)\xrarrow{\e}\ZZ\to0
\]
given as follows.

Let $F_k(\Gam_0)$ be the free abelian group generated by all
isomorphism classes of codimension $k$ objects in $\Fat$ over
$\Gam_0$ together with an orientation. In other words, generators
of $F_k(\Gam_0)$ are given by morphisms
\[
    f:\Gam\to\Gam_0
\]
where $\Gam$ is a ribbon graph of codimension $k$ together with an
orientation on $\Gam$. Two such objects $f_i:\Gam_i\to\Gam_0$ for
$i=1,2$ are \emph{isomorphic} if there is an orientation
preserving isomorphism $g:\Gam_1\to\Gam_2$ so that $f_2\circ
g=f_1$. As usual, we equate reversal of orientation with reversal
of sign. In particular, $F_n(\Gam_0)$ has rank $1$ with two
generators corresponding to the two possible orientations of
$\Gam_0$.

The boundary map $d:F_k(\Gam_0)\to F_{k-1}(\Gam_0)$ is given by
\[
    d[f:\Gam\to\Gam_0]=\sum [f\circ g_i:\Gam_i\to\Gam_0]
\]
where the sum is taken over all right equivalence classes of
morphisms
\[
    g_i:\Gam_i\to \Gam
\]
which collapse only one edge. We take the unique orientation on
each $\Gam_i$ which induces the given orientation on $\Gam$.

The augmentation map $\e:F_0(\Gam_0)\to\ZZ$ is given by
\[\e[\Gam\to\Gam_0]=o(\Gam)=\pm1.\]

\begin{prop}\label{prop:forested graph complex is acyclic} Suppose
that $\Gam_0$ is trivalent except for $r$ vertices
$v_1,\cdots,v_r$ which have codimensions $n_1,\cdots,n_r$, resp.
Then $F_\ast(\Gam_0)$ is based chain isomorphic to the tensor
product
\[
    F_\ast(\Gam_0)\cong C_\ast(K^{n_1})\otimes
    C_\ast(K^{n_2})\otimes\cdots\otimes C_\ast(K^{n_r})
\]
where $C_\ast(K^m)$ is the cellular chain complex of the $m$
dimensional Stasheff polyhedron $K^m$. In particular,
$F_\ast(\Gam_0)$ is acyclic.
\end{prop}

Proposition \ref{prop:forested graph complex is acyclic} follow
from well-known properties of the Stasheff polyhedron which we
will review shortly. Suppose for the moment that this is true.

For each $\Gam_0$ there is a natural augmented chain map
$p:F_\ast(\Gam_0)\to \G_\ast^\ZZ$ given by
\[
    p[f:\Gam\to\Gam_0]=\<\Gam\>.
\]
A morphism $g:\Gam_0\to\Gam_1$ induces a chain map
$g_\ast:F_\ast(\Gam_0)\to F_\ast(\Gam_1)$ by
\[
    g_\ast[f:\Gam\to\Gam_0]=[g\circ f:\Gam\to\Gam_1].
\]
This is a chain map over $\G_\ast^\ZZ$ in the sense that $p\circ
g_\ast=g_\ast$. Therefore, $F_\ast$ is a functor from $\Fat$ to
the category of acyclic augmented chain complexes over
$\G_\ast^\ZZ$. In other words, it is an acyclic carrier. We call
$F_\ast$ the \emph{forest carrier}.

The acyclic carrier $F_\ast$ carries a unique (up to homotopy)
chain map
\begin{equation}\label{eq:chain map given by F}
    \f_\ast:C_\ast(\Fat)\to \G_\ast^\ZZ
\end{equation}
where $C_\ast(\Fat)$ is the cellular chain complex of the category
of ribbon graphs. (So $C_n(\Fat)$ is the free abelian group
generated by all $n$ simplices
\[
    \Gam_\ast=(\Gam_0\to\Gam_1\to\cdots\to\Gam_n)
\]
in the nerve of $\Fat$.)

\begin{thm}\label{thm:F carries a rational equivalence}
Any chain map (\ref{eq:chain map given by F}) carried by $F_\ast$
is a rational homotopy equivalence.
\end{thm}

\begin{rem}
If we consider the forest carrier $F_\ast$ as a diagram of chain
complexes over $\G_\ast^\ZZ$ we see that there is an induced chain
map from the homotopy pushout of this diagram into $\G_\ast^\ZZ$.
This homotopy pushout is the \emph{forested graph complex}
\[
    C_\ast(\Fat;F_\ast)=\bigoplus_n\bigoplus_{\Gam_\ast\in\N_n\Fat}
    \s^n\ZZ(\Gam_\ast)\otimes F_\ast(\Gam_0)
\]
where $\ZZ(\Gam_\ast)$ is the free abelian group of rank one
generated by $(\Gam_\ast)$ and $\s^n$ is the $n$-fold suspension
operator. Since $F_\ast$ is acyclic, $C_\ast(\Fat;F_\ast)\simeq
C_\ast(\Fat;\ZZ)$. This gives the following diagram.
\[
    C_\ast(\Fat;\ZZ)\xleftarrow{\simeq}C_\ast(\Fat;F_\ast)\xrarrow{p}
    \G_\ast^\ZZ.
\]
We are claiming that the right hand arrow, given by
$p:F_\ast\to\G_\ast^\ZZ$ for $n=0$ and zero for $n>0$, is a
rational homotopy equivalence.
\end{rem}

\begin{thm}\label{thm:W ll pulls back to integer cohomology class}
The dual Kontsevich cycles pull back to well-defined integer
cohomology classes
\[
    \f^\ast[W_\ll^\ast]\in H^{2|\ll|}(\Fat;\ZZ)\cong \prod_{g,s}
    H^{2|\ll|}(M_g^s;\ZZ).
\]
\end{thm}

To prove Theorem \ref{thm:F carries a rational equivalence} we
will obtain a description of the chain map (\ref{eq:chain map
given by F}) and construct an explicit rational homotopy inverse
using ``dual cells.'' The proof assumes Proposition
\ref{prop:forested graph complex is acyclic}.

\subsection{Dual cells}

``Dual cells'' are elements of $C_\ast(\Fat)$ associated to every
generator $\<\Gam\>\in\G_n^\ZZ$. Every oriented ribbon graph
$\Gam$ has many dual cells but we will see that each of them is
necessarily mapped to $(-1)^{\binom{n+1}{2}}\<\Gam\>$ by any chain
map carried by the forest carrier.

A rational inverse
\[
    \psi:\G_\ast^\ZZ\otimes\QQ\to C_\ast(\Fat;\QQ)
\]
is given by mapping each rational generator $\<\Gam\>$ to the
average dual cell in a finite model for $C_\ast(\Fat;\QQ)$ given
by choosing one object from every isomorphism class of ribbon
graphs and taking the average over all possible dual cells which
lie in this finite model.

The composition $\f\psi$ will be the identity mapping on
$\G_\ast^\ZZ\otimes\QQ=\G_\ast\otimes\QQ$ and the composition
$\f\psi$ will be homotopic to the identity on the finite model
since it is carried by the ``identity carrier.'' The identity
carrier is a canonical acyclic carrier which carries the identity
map on the cellular chain complex of any small category. (See
Lemma \ref{lem:thin to G to thin}.)

Suppose that $\Gam$ is an oriented ribbon graph of codimension
$n$. Then a \emph{dual cell} for $\Gam$ is given by choosing one
representative from every isomorphism class of ribbon graphs over
$\Gam$ (taking the identity map on $\Gam$ as one representative).
Consider all $n$ simplices
\[
    \Gam_\ast=(\Gam_0\to\cdots\to\Gam_n=\Gam)
\]
where $\Gam_i$ is a representative of codimension $i$. Then the
dual cell is given by the signed sum of all of these
$n$-simplices:
\[
    D(\Gam)=\sum o(\Gam_\ast)(\Gam_\ast)\in C_n(\Fat)
\]
where the sign $o(\Gam_\ast)=\pm1$ is positive iff the given
orientation of $\Gam$ agrees with the one induced from the natural
orientation of the trivalent graph $\Gam_0$. Note that, if we
reverse the orientation of $\Gam$, this sign will change. So,
\[
    D(-\Gam)=-D(\Gam).
\]
It is also trivial to see that, in the case $n=0$, we have
$D(\Gam)=\Gam$ assuming that we take the natural orientation on
$\Gam$.

\begin{lem}\label{lem:boundary of dual cell}
The boundary of the dual cell is, up to sign, a sum of dual cells:
\[
    dD(\Gam)=(-1)^{n}\sum D(\Gam')
\]
where the sum is taken over all chosen
representatives\[\Gam'\to\Gam\]of all isomorphism classes of
ribbon graphs of codimension $n-1$ over $\Gam$. We take the
orientation on $\Gam'$ induced by the given map.
\end{lem}

\begin{proof}
By Proposition \ref{prop:left cancellation}, the set of
isomorphism classes of objects over each $\Gam'$ maps
monomorphically into the set of isomorphism classes of objects
over $\Gam$. Therefore, the given choices of representatives for
$\Gam$ gives a complete set of representatives for the objects
over $\Gam'$ and $D(\Gam')$ is defined.

The boundary of $D(\Gam)$ is given by
\begin{multline*}
    dD(\Gam)=\sum_{\Gam_\ast}\sum_{i=0}^{n-1}
    (-1)^io(\Gam_\ast)(\Gam_0,\cdots,\widehat{\Gam}_i,\cdots,\Gam_n=\Gam)\\
    +(-1)^n\sum_{\Gam_\ast}o(\Gam_\ast)(\Gam_0,\cdots,\Gam_{n-1})
\end{multline*}
However, the double sum is zero since, when $\Gam_i$ is deleted,
there are exactly two ways to fill in the blank and these give
opposite signs for $o(\Gam_\ast)$. The second sum is equal to the
sum of $D(\Gam_{n-1})$ for all possible $\Gam_{n-1}$.
\end{proof}

\begin{lem}\label{lem:D(Gam) is mapped to pm Gam}
Given any oriented ribbon graph $\Gam$ of codimension $n$, any
dual cell $D(\Gam)\in C_\ast(\Fat)$ and any augmented chain map
$\f:C_\ast(\Fat)\to\G_\ast^\ZZ$ carried by the forest carrier
$F_\ast$ we will have \[
    \f(D(\Gam))=(-1)^{\binom{n+1}{2}}
    \<\Gam\>.
    \]
\end{lem}

\begin{proof}
This will be by induction on $n$. Suppose that $n=0$. Then
\[
    \f(D(\Gam))=\f(\Gam)=\<\Gam\>
\]
since the identity map $[\Gam\to\Gam]$ is the unique element of
$F_0(\Gam)$ with augmentation equal to $1$.

Now suppose the statement holds for $n-1$. Then by Lemma
\ref{lem:D(Gam) is mapped to pm Gam} we have
\[
    \f dD(\Gam)=(-1)^n\f\sum
    D(\Gam')=(-1)^{n+\binom{n}{2}}\sum\<\Gam'\>=(-1)^{\binom{n+1}{2}}d\<\Gam\>.
\]
However, $F_{n+1}(\Gam)=0$. So, the value of $\f D(\Gam)$ in
$F_n(\Gam)$ is uniquely determined by the value of its boundary.
This forces $\f D(\Gam)$ to be $(-1)^{\binom{n+1}{2}}\<\Gam\>$.
\end{proof}

To construct a rational inverse for the chain map $\f$ we choose a
finite model for $\Fat$. Let $\thin$ be a full subcategory of
$\Fat$ that contains exactly one object from every isomorphism
class. Then $\thin$ is a deformation retract of $\Fat$ and the
cellular chain complex of $\thin$ is a deformation retract of
$C_\ast(\Fat)$.

If $\Gam$ is any oriented ribbon graph of codimension $n$, let
$\overline{D}(\Gam)\in C_n(\thin;\QQ)$ be the \emph{average dual
cell} of $\Gam$ given by
\begin{equation}\label{eq:average dual cell}
    \overline{D}(\Gam)=\sum
    \frac{o(\Gam_\ast)}{|\Aut(\Gam_0,\cdots,\Gam_{n})|}
    (\Gam_0\to\cdots\to\Gam_n\xrarrow{\approx}\Gam)
\end{equation}
where
\[
    \Aut(\Gam_0,\cdots,\Gam_{n})
    =\Aut(\Gam_0)\times\cdots\times\Aut(\Gam_{n})
\]
and the sum is taken over all possible sequences of morphism in
$\thin/\Gam$ so that $\Gam_i$ has codimension $i$ for each $i$.

\begin{lem}\label{lem:d ave dual cell} If $d\<\Gam\>=\sum \<\Gam'\>$
in $\G_\ast^\ZZ$ then $d\overline{D}(\Gam)=(-1)^n\sum
\overline{D}(\Gam')$.
\end{lem}

\begin{proof}
The proof is analogous to the proof of Lemma \ref{lem:boundary of
dual cell}. We just need to realize that when $\Gam_i$ is deleted,
there are $|\Aut(\Gam_i)|$ ways to put it back in as an isomorphic
copy. Consequently,
\begin{multline*}
    d\overline{D}(\Gam)=\sum_{i=0}^{n-1} (-1)^i\sum
    \frac{o(\Gam_\ast)}{|\Aut(\Gam_0,\cdots,\widehat{\Gam}_i,\cdots,\Gam_{n})|}
    (\Gam_0\to\cdots\widehat{\Gam}_i\cdots\to\Gam_n\xrarrow{\approx}\Gam)\\
    +(-1)^n\sum
    \frac{o(\Gam_\ast)}{|\Aut(\Gam_0,\cdots,\Gam_{n-1})|}
    (\Gam_0\to\cdots\to\Gam_{n-1}\to\Gam).
\end{multline*}
Then, in the second sum, the morphism $\Gam_{n-1}\to\Gam$ can be
uniquely factored through some $\Gam'$ making it into a sum of
terms of the form $(-1)^n \overline{D}(\Gam')$.
\end{proof}

This lemma says that we have a chain map
\[
    \psi:\G_\ast^\ZZ\otimes\QQ\to C_\ast(\thin;\QQ)
\]
given by $\psi\<\Gam\>=(-1)^{\binom{n+1}{2}}\overline{D}(\Gam)$.
Lemma \ref{lem:D(Gam) is mapped to pm Gam} gives us:

\begin{thm}\label{thm:G to thin to Fat to G is identity} The composition
\[
    \G_\ast\otimes\QQ\xrarrow{\psi} C_\ast(\thin;\QQ)\hookrightarrow
    C_\ast(\Fat;\QQ)\xrarrow{\f}\G_\ast\otimes\QQ
\]
is the identity map on $\G_\ast\otimes\QQ=\G_\ast^\ZZ\otimes\QQ$
for any chain map $\f$ carried by the forest carrier.
\end{thm}

Finally, Theorem \ref{thm:F carries a rational equivalence}
follows from the following.

\begin{lem}\label{lem:thin to G to thin} The composition
\[
    C_\ast(\thin;\QQ)\hookrightarrow
    C_\ast(\Fat;\QQ)\xrarrow{\f}\G_\ast\otimes\QQ\xrarrow{\psi} C_\ast(\thin;\QQ)
\]
is homotopic to the identity map.
\end{lem}

\begin{proof}
To show that two chain maps are homotopic it suffices to construct
an acyclic carrier that carries both of them. In this case it will
be the ``identity carrier.''

Let $X$ be any object of any small category $\A$. Then the
category $\A/X$ of objects over $X$ is contractible since it has a
terminal object $id:X\to X$ and any morphism $X\to Y$ induces a
functor $\A/X\to\A/Y$. Consequently, the cellular chain complex
$C_\ast(\A/X)$ is an acyclic carrier from $C_\ast(\A)$ to itself,
i.e., a functor from $\A$ into the the category of augmented
acyclic chain complexes over $C_\ast(\A)$. We call this the
\emph{identity carrier} since it carries the identity morphism on
$C_\ast(\A)$.

In order to show that $\psi\circ \f$ is carried by the identity
carrier we need to show that the chain map $\psi$ is covered by a
mapping from the forest carrier to the identity carrier, i.e., we
need a natural commuting diagram as follows for all objects $\Gam$
in $\thin$.
\[
    \begin{diagram}
\node{F_\ast(\Gam)\otimes\QQ}\arrow{e,t}{\widetilde{\psi}}\arrow{s,l}{p}
\node{C_\ast(\thin/\Gam;\QQ)}\arrow{s,r}{q}\\
\node{\G_\ast\otimes\QQ}\arrow{e,t}{\psi}\node{C_\ast(\thin;\QQ)}
    \end{diagram}
\]
Recall that the generators of $F_\ast(\Gam)$ are isomorphism
classes $[f:\Gam'/\to\Gam]$ of oriented ribbon graphs over $\Gam$.
The projection $p$ sends this to $\<\Gam'\>\in \G_\ast\otimes\QQ$
which then goes to the average dual cell $\overline{D}(\Gam')\in
C_\ast(\thin;\QQ)$. But the morphism $f:\Gam'\to\Gam$ makes all
the terms in the definition (\ref{eq:average dual cell}) of
$\overline{D}(\Gam')$ into simplices in $\thin/\Gam$ and therefore
defines a lifting
\[
    \widetilde{\psi}[f:\Gam'\to\Gam]=f_\ast(\overline{D}(\Gam'))
    \] of $\overline{D}(\Gam')$ to $C_\ast(\thin/\Gam;\QQ)$ as
    required.
\end{proof}

 %\vfill\eject

\setcounter{section}{1}

\section{Stasheff associahedra}

We use several different versions of the Stasheff associahedron:
the convex $n$ dimensional polyhedron $K^n$, the category
$\A_n=\simp K^{n-3}$ of simplices in $K^{n-3}$ and the simplicial
nerve of $\A_n$ which is a triangulation of the polyhedron
$K^{n-3}$.

An outline of this section:
\begin{enumerate}
    \item Stasheff polyhedron $K^n$.
    \item The category $\A_{n+3}$.
    \item Orientation of $K^n$.
    \item Orientation of $K^{odd}$.
    \item Proof of Proposition \ref{prop:forested graph complex is
acyclic}.
\end{enumerate}

\subsection{Stasheff polyhedron $K^n$}

\begin{thm}[Stasheff \cite{[Stasheff63:Polyhedron]}]\label{thm:existence of Stasheff polyhedron}
There is an $n$-dimensional convex polyhedron $K^n$ whose elements
correspond to isomorphism classes of planar metric trees with
$n+3$ leaves of fixed length and up to $n$ internal edges of
variable length $\leq1$. Two planar metric trees lie in the same
open face of $K^n$ if and only if, after collapsing all internal
edges of length $<1$, they become isomorphic fixing the leaves.
\end{thm}

Let $C_\ast(K^n)$ be the cellular chain complex of $K^n$. Then
$C_0(K^n)$ is freely generated by isomorphism classes of trivalent
planar trees with $n+3$ given leaves. Recall that any such tree
has a natural orientation. More generally we have the following.

\begin{prop}\label{prop:cellular chain complex of Kn} For all
$0\leq k\leq n$, $C_k(K^n)$ is generated by isomorphism classes
$[T]$ of oriented planar trees $T$ with $n+3$ fixed leaves and
$n-k$ internal edges. The boundary map is given by
\[
    d[T]=\sum_{[T',e]} [T']
\]
where the sum is taken over all isomorphism classes of pairs
$(T',e)$ where $e$ is an edge in an oriented tree $T'$ so that
$T\cong T'/e$ with the induced orientation.
\end{prop}

\begin{rem}\label{rem:geometric orientation is chosen}
This is more or less a tautology since we choose the geometric
orientation of the faces to make this algebraic statement true. We
also note that trees with fixed leaves have no nontrivial
automorphisms. Therefore,
\[
    \<T\>=|\Aut(T)|[T]=[T].
\]
\end{rem}

\begin{proof}
Let $T_0$ be a trivalent planar tree with $n+3$ fixed leaves
$h_0,\cdots,h_{n+2}$ (in cyclic order) and $n$ internal edges
$e_1,\cdots,e_n$. Since $T_0$ is trivalent, it has an intrinsic
orientation. If we collapse the edges $e_1,\cdots,e_k$ in that
order we get a tree $T_k$ of codimension $k$ with the induced
orientation. The trees $T_k$ for various $k$ are related by
\[
    T_k=T_{k-1}/e_k
\]
The $k$ dimensional face of $K^n$ corresponding to the tree $T_k$
consists of isomorphism class of trees $T$ having edges
$e_1,\cdots,e_k$ of variable length $<1$ and the other edges of
length equal to $1$. Let $x_1,\cdots,x_k$ be the lengths of
$e_1,\cdots,e_k$ then we choose the geometric orientation of this
face of $K^n$ (in a neighborhood of the vertex $T_0$) by taking
these coordinated in opposite order:
\[
    (x_k,\cdots,x_1).
\]
When $x_k$ reaches $1$ we get to the face corresponding to
$T_{k-1}$. The orientation of the $k-1$ face is therefore related
to that of the $k$-face by the ``first vector points outward''
rule which is standard.
\end{proof}

\subsection{The category $\A_n$}

Let $\A_{n+3}$ be the category of faces of $K^n$ with inclusion
maps as morphisms. Then the geometric realization of the nerve of
$\A_{n+3}$ is homeomorphic to $K^n$. A homeomorphism
\[
    \f:|\A_{n+3}|\to K^n
\]
is given by sending each object (=vertex in the nerve) to some
point in the interior of the corresponding face of $K_n$ and
extending linearly over the simplices.

Take an $n$-simplex
\begin{equation}\label{eq:maximal dim simplex in Kn}
   T_\ast=(T_0\to T_1\to \cdots\to T_n)
\end{equation}
in $\A_{n+3}$ which is \emph{nondegenerate} in the sense that
$T_0$ is trivalent and each $T_i$ is obtained from $T_{i-1}$ by
collapsing one edge. Then, by induction, each tree $T_i$ obtains
an induced orientation from the intrinsic orientation of $T_0$.
The orientation of the tree $T_n$ gives a geometric orientation of
the polyhedron $K^n$ as explained in the proof of Proposition
\ref{prop:cellular chain complex of Kn}.

Suppose that an orientation for the tree $T_n$ representing the
top cell of $K^n$ is given. (For example, when $n$ is even, we can
take the intrinsic orientation of $T_n$.) Then the \emph{algebraic
orientation} $o(T_\ast)$ of the $n$-simplex $T_\ast$ in
(\ref{eq:maximal dim simplex in Kn}) is defined to be $\pm1$
depending on whether the given orientation of $T_n$ is equal to
the one induced from $T_0$. The \emph{dual cell} $D(T_n)$ is
defined by
\[
    D(T_n)=\sum o(T_\ast)T_\ast\in C_n(\A_{n+3})
\]
where the sum is over all nondegenerate $n$-simplices $T_\ast$.
This is an $n$-chain in the cellular chain complex of $\A_{n+3}$.
It represents a triangulation of the top cell of $K^n$ with
algebraic orientation which does not agree with the geometric
orientation.

\subsection{Orientation of $K^n$}

\begin{lem}\label{lem:simplicial vs geometric orientation}
The embedding
\[
    \Delta^n\to K^n
\]
given by the $n$ simplex (\ref{eq:maximal dim simplex in Kn}) has
degree $(-1)^{\binom{n+1}{2}}$ with respect to the orientation of
$K^n$ corresponding to the orientation of $T_n$ induced from that
of $T_0$.
\end{lem}

\begin{rem}\label{rem:orientation of K agrees with psi}
This implies that, given any orientation of $T_n$, the
corresponding geometric orientation of $K^n$ agrees with
\[
    \psi[T_n]:=(-1)^{\binom{n+1}{2}}{D}(T_n).
\]
\end{rem}

\begin{proof} For any $0\leq k\leq n$ we claim that the map
\[
    \s^k:\Delta^k\to F_k
\]
given by $T_0\to\cdots\to T_k$ has degree $(-1)^{\binom{k+1}{2}}$
where $F_k$ is the face of $K^n$ corresponding to $T_k$ with the
induced orientation.

This statement hold for $k=0$. Suppose it holds for $k-1$. Since
the front $k-1$ face $\Delta^{k-1}$ is opposite the $k$-th vertex
$T_k$, its orientation is equal to $(-1)^k$ times the induced
orientation from $\Delta^k$. Whereas, we are orienting the faces
to make Proposition \ref{prop:cellular chain complex of Kn} true,
i.e., the orientation of $F_{k-1}$ is the one induced from $F_k$.
Consequently the degree of $\s^k$ is
\[
    (-1)^{\binom{k}{2}}(-1)^k=(-1)^{\binom{k+1}{2}}.
\]
Putting $k=n$ we get the lemma.
\end{proof}

Now suppose that $n=2k$. Then $T_{2k}$ has an intrinsic
orientation which determines an intrinsic orientation for the
polyhedron $K^{2k}$.

\begin{thm}\label{thm:orientation of simplices in Keven} The sign
of the embedding
\[
    \s^{2k}:\Delta^{2k}\to K^{2k}
\]
given by the $2k$ simplex $(T_0\to T_1\to\cdots \to T_{2k})$ is as
given in \cite{[I:MMM_and_Witten]}.
\end{thm}

In \cite{[I:MMM_and_Witten]} the sign of $\s^{2k}$ is defined to
be the sign of the permutation of $2k+3$ letters:
\[
    \sgn(\s^n)=\sgn(a_1,a_2,a_3,b_1,b_2,\cdots,b_{2k})
\]
where $a_i,b_j$ are the \emph{regions} (components of the
complement of the tree in the disk in which it can be embedded
with leaves on the boundary) given as follows. The regions which
bound the edge $e_1$ (so that $T_1=T_0/e_1$) are $a_1,a_3$ (in
either order). Let the two other regions which touch $e_1$ be
$a_2,b_1$ so that they are $a_1,a_2,a_3,b_1$ in cyclic order. For
$i\geq 2$ let $v_i$ be the vertex of $e_i$ furthest away from
$e_1$ and let $b_i$ be the region which touches $e_i$ at the point
$v_i$. Then $\sgn(\s^n)$ is equal to the sign of the permutation
of $a_1,a_2,a_3,b_1,\cdots,b_{2k}$ which puts these regions into
the correct cyclic order

\begin{proof}
We have three definitions of the orientation of $K^{2k}$.
\begin{enumerate}
    \item[(a)] The intrinsic orientation of $K^{2k}$ induced from the intrinsic orientation of $T_{2k}$.
    \item[(b)] The orientation of $K^{2k}$ induced from that of $T_0$
    by the sequence of maps $T_0\to\cdots\to T_{2k}$.
    \item[(c)] The orientation of $K^{2k}$ induced from $\Delta^{2k}$
    by the map $\s^{2k}$.
\end{enumerate}
The statement we are trying to prove is that (a) and (c) differ by
the sign convention given above. We know, by Lemma
\ref{lem:simplicial vs geometric orientation}, that (b) and (c)
differ by the sign
\[
    (-1)^{\binom{2k+1}{2}}=(-1)^k.
\]
Therefore, it suffices to show that the sign difference between
(a) and (b), which is equal to $o(T_\ast)$ by definition, is given
by
\begin{equation}\label{eq:reversed sign of simplex}
   o(T_\ast)= (-1)^k\sgn(\s^{2k})=\sgn(a_1,a_2,a_3,b_{2k},\cdots,b_1).
\end{equation}
But this is a special case of the following lemma.
\end{proof}

\begin{lem}\label{lem:relative sign of collapses tree, general case}
Suppose that $T_0$ is a planar tree with internal edges
$e_1,\cdots,e_{2k}$ and $2k+1$ vertices all trivalent except for
$v_0$ which has valence $2n+1$. Let $v_i$ be the vertex on $e_i$
furthest away from $v_0$. Let $T_1,\cdots,T_{2k}$ be given by
$T_i=T_{i-1}/e_i$. Then the difference $o(T_\ast)$ between the
intrinsic orientation of $T_{2k}$ and the one induced from $T_0$
is given by
\begin{equation}\label{eq:sign convention for back face}
    o(T_\ast)=(-1)^k\sgn(a_1,\cdots,a_{2n+3},b_1,\cdots,b_{2k})
\end{equation}
where $a_i$ are the regions around $v_0$ and $b_i$ is the region
which touched $e_i$ only at $v_i$.
\end{lem}

\begin{proof}
First we claim that the truth value of this statement remains
unchanged if we permute the edges $e_1,\cdots,e_{2k}$. To see this
suppose we change the order of $e_i,e_{i+1}$. Then $T_i$ will
become a different tree and the orientations of
$T_{i+1},\cdots,T_{2k}$ will be reversed. This comes from the
proof that $\d^2=0$ in graph homology. Switching $e_i,e_{i+1}$
will also transpose the labels $b_i,b_{i+1}$. So the sign
(\ref{eq:sign convention for back face}) will also change so the
relative sign remains unchanged.

By permuting the order of the edges we may assume that each tree
$T_i$ has only one vertex $v_0$ which is not trivalent and
$e_{i+1}$ becomes an edge in $T_i$ connecting $v_0$ to $v_{i+1}$.
Consequently, $T_{2k-2}$ has an intrinsic orientation and, by
induction on $k$, this orientation differs from the one induced
from $T_0$ by
\begin{equation}\label{eq:sign convention for middle face}
    o(T_0,\cdots,T_{2k-2})=(-1)^{k-1}\sgn(a_1,\cdots,a_{2n+3},b_1,\cdots,b_{2k-2})
\end{equation}

To prove that (\ref{eq:sign convention for middle face}) implies
(\ref{eq:sign convention for back face}) we look at the difference
between the two signs and the difference between the intrinsic
orientation of $T_{2k}$ and the one induced from the intrinsic
orientation of $T_{2k-2}$. The statement that we need is exactly
the statement of the lemma in the case $k=1$.

Now we assume that $k=1$. As before, we may assume that $e_1$
connects $v_0$ to $v_1$ but there are two possibilities for $e_2$.

\underline{Case 1.} $e_2$ connects $v_0$ and $v_2$.
%*
\begin{figure}
\includegraphics{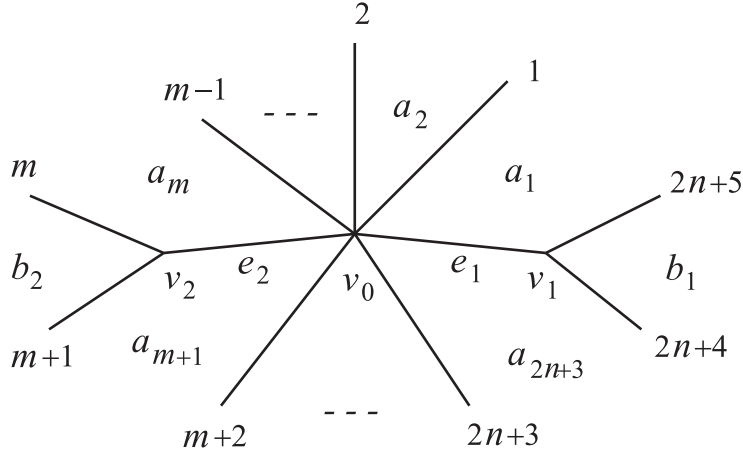}
\caption{Case 1. $\sgn(a_1,\cdots,a_{2n+3},b_1,b_2)=(-1)^m$ where
$1\leq m\leq 2n+2$.}\label{fig:tree01}%\treeone
\end{figure}
%*

If we let $a_1,\cdots,a_m$ be the regions at $v_0$ from $e_1$ to
$e_2$ as shown in Figure \ref{fig:tree01} then the intrinsic
orientation of $T_0$ and induced orientations of $T_1,T_2$
are
\begin{enumerate}
    \item[($T_0$)] ${\bf v_0e_1^-}h_1h_2\cdots h_{m-1}e_2^-h_{m+2}\cdots
    h_{2n+3}\ {\bf v_1e_1^+}h_{2n+4}h_{2n+5}\ v_2e_2^+h_mh_{m+1}$
    \item[($T_1$)] $-{\bf v_0}h_1h_2\cdots h_{m-1}{\bf e_2^-}h_{m+2}\cdots
    h_{2n+5}\ {\bf v_2e_2^+}h_mh_{m+1}$
    \[
        =(-1)^m{\bf v_0e_2^-v_2e_2^+}h_1\cdots h_{2n+5}
    \]
    \item[($T_2$)] $(-1)^{m-1}v_0h_1h_2\cdots
    h_{2n+5}$
\end{enumerate}
where the leaves are labelled $h_1,\cdots,h_{2n+5}$ in
counterclockwise order and $e_i^-,e_i^+$ are the halves of $e_i$
closer/further from $v_0$. We see that the induced orientation of
$T_1$ differs from the intrinsic orientation by $(-1)^{m-1}$.
However, the regions $a_1,\cdots,a_{2n+3},b_1,b_2$ are arranged in
the cyclic order
\[
    [b_1,a_1,\cdots,a_m,b_2,a_{m+1},\cdots,a_{2n+3}]
\]
so the permutation sign is
\[
    \sgn(a_1,\cdots,a_{2n+3},b_1,b_2)=(-1)^m
\]
Multiplying by $(-1)^k=-1$ we get (\ref{eq:sign convention for
back face}).

\underline{Case 2}. $e_2$ connects $v_1$ and $v_2$.
%*
\begin{figure}
\includegraphics{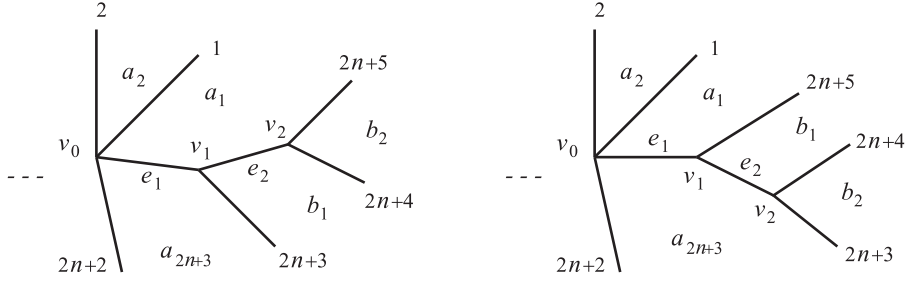}
\caption{Case 2a: $b_1b_2$ (on left) and Case 2b: $b_2b_1$ (on
right).}\label{fig:tree02}%\treetwo
\end{figure}
%*

Here there are two subcases as shown in Figure
\ref{fig:tree02}.\begin{enumerate}
    \item[2a] $b_1$ is clockwise from $b_2$.
    \item[2b] $b_2$ is clockwise from $b_1$.
\end{enumerate}
In subcase 2a we have the following orientations on $T_0,T_1,T_2$
induced from the intrinsic orientation of $T_0$.
\begin{enumerate}
    \item[($T_0$)] ${\bf v_0e_1^-}h_1h_2\cdots h_{2n+2}
    \ {\bf v_1e_1^+}h_{2n+3}e_2^-\ v_2e_2^+h_{2n+4}h_{2n+5}$
    \item[($T_1$)] $-{\bf v_0}h_1h_2\cdots
    h_{2n+3}\ {\bf e_2^-v_2e_2^+}h_{2n+4}h_{2n+5}={\bf v_0e_2^-v_2e_2^+}h_1\cdots h_{2n+5}$
    \item[($T_2$)] $-v_0h_1h_2\cdots
    h_{2n+5}$
\end{enumerate}
The induced orientation of $T_2$ is negative the natural
orientation. But
\[
    \sgn(a_1,\cdots,a_{2n+3},b_1,b_2)=+1
\]
since the regions are in the correct cyclic order. Thus the lemma
holds in this case.

In subcase 2b the induced orientations on $T_0,T_1,T_2$ are:
\begin{enumerate}
    \item[($T_0$)] ${\bf v_0e_1^-}h_1h_2\cdots h_{2n+2}
    \ {\bf v_1e_1^+}e_2^-h_{2n+5}\ v_2e_2^+h_{2n+3}h_{2n+4}$
    \item[($T_1$)] $-{\bf v_0}h_1h_2\cdots
    h_{2n+2}{\bf e_2^-}h_{2n+5}\ {\bf v_2e_2^+}h_{2n+3}h_{2n+4}=-{\bf v_0e_2^-v_2e_2^+}h_1\cdots h_{2n+5}$
    \item[($T_2$)] $v_0h_1h_2\cdots
    h_{2n+5}$
\end{enumerate}
The induced orientation on $T_2$ is equal to the natural
orientation. But
\[
    \sgn(a_1,\cdots,a_{2n+3},b_1,b_2)=-1
\]
so the lemma holds in this final case.
\end{proof}

\subsection{Orientation of $K^{odd}$}

Suppose that $n$ is odd. Then the planar trees in $K^n$ have an
even number ($n+3$) of leaves. An orientation of these trees is
given by taking a fixed cyclic ordering of these leaves, starting
at one point and going counterclockwise. This gives an orientation
of a tree with one vertex and thus of the top cell of $K^n$. For
any other tree $T\in K^n$ an orientation is given by choosing an
ordering of the internal edges of $T$. The orientation given by
this ordering is the one which induces the chosen orientation on
the one vertex tree by collapsing the edges in order.

We can now prove Proposition \ref{prop:forested graph complex is
acyclic}.

\subsection{Proof of Proposition \ref{prop:forested graph complex is
acyclic}}

Recall that we have a ribbon graph $\Gam_0$ all of whose vertices
are trivalent except for $v_1,\cdots,v_r$ which have valence
$n_1+3,\cdots,n_r+3$. For each $i$ we choose a cyclic ordering of
the half-edges at $v_i$ and take the orientation on $\Gam_0$ given
by $v_1$ followed by its half-edges, $v_2$ followed by its
half-edges, etc.

Not suppose that for each $i$, $T_i$ is a planar tree with $n_i+3$
leaves representing an $m_i$-face of $K^{n_i}$. Thus $T_i$ has
$n_i-m_i$ internal edges. Choose an ordering of these edges. Let
$\Gam_0(T_1,\cdots,T_r)$ denote the graph obtained from $\Gam_0$
by replacing $v_i$ by $T_i$. Then the $T_i$ will form a forest in
$\Gam_0$. Let $e_{i1},e_{i2},\cdots$ be the internal edges of
$T_i$. Take the orientation on $\Gam_0(T_1,\cdots,T_r)$ so that,
if the edges $e_{ij}$ are collapsed in lexicographic order, we get
the chosen orientation on $\Gam_0$.

Let
\[
    \f:C_{m_1}(K^{n_1})\otimes\cdots\otimes C_{m_r}(K^{n_r})\to
    F_m(\Gam_0)
\]
where $m=m_1+\cdots +m_r$ be given by
\[
    \f([T_1]\otimes \cdots\otimes [T_r])=[\Gam_0(T_1,\cdots,T_r)\to\Gam_0].
\]
We claim that this gives a chain isomorphism
\[
    \f:C_\ast(K^{n_1})\otimes\cdots\otimes C_\ast(K^{n_r})\to
    F_\ast(\Gam_0).
\]
Since $\f$ sends basis elements to basis elements, it suffices to
show that $\f$ is a chain map. But this is straightforward.

The boundary of $[T_1]\otimes\cdots\otimes [T_r]$ is, by
Proposition \ref{prop:cellular chain complex of Kn}, equal to
\[
    \sum_{i=1}^r(-1)^{m_1+\cdots m_{i-1}}\sum [T_1]\otimes\cdots\otimes [T_{i-1}]\otimes
    [T_i']\otimes[T_{i+1}]\otimes\cdots\otimes [T_r]
\]
where the second sum is over all pairs $(T_i',e_{i0})$ so that
$T_i'/e_{i0}\cong T_i$. But $(-1)^{m_1+\cdots m_{i-1}}$ is also
the sign of the permutation which brings the edge $e_{i0}$ to the
beginning in the ordering of all edges of
$\Gam_0(T_1,\cdots,T_i',\cdots,T_r)$. So
\[
    d\f([T_1]\otimes\cdots\otimes[T_r])=\sum (-1)^{m_1+\cdots
    m_{i-1}}\f([T_1]\cdots[T_i']\cdots[T_r])
\]
as required.

 %\vfill\eject

\setcounter{section}{2}

\section{Miller-Morita-Mumford classes} % sec 03

We review the definition of the adjusted Miller-Morita-Mumford
classes. We pull them back to the graph cohomology complex
$\G_\ast$ and we describe what happens when we evaluate them on
the dual Kontsevich cycles.

\begin{enumerate}
    \item Cyclic set cocycle.
    \item Adjusted MMM classes in $H^{2k}(\G_\ast;\QQ)$.
    \item Cup products of adjusted MMM classes.
    \item Computing the numbers $b_{n_\ast}^{k_\ast}$.
    \item Kontsevich cycles in terms of MMM classes.
    \item Computing $a_\ll^\mu$.
\end{enumerate}

Suppose that \[\Sig_g^s\to E\xrarrow{p} B\] is a compact manifold
bundle where $\Sig_g^s$ is an oriented connected surface of genus
$g$ and $s$ unordered distinguished points. Let
\[
    \pi:\widetilde{B}\to B
\]
be the $s$-fold covering space given by the $s$ distinguished
point in each fiber of $E$. Then the vertical tangent bundle of
$E$ is an oriented $2$-plane bundle and therefore has an Euler
class $e(E)\in H^2(E;\ZZ)$. The push-down of the $k+1$-st power of
this class is the \emph{Miller-Morita-Mumford class}
\[
    \k_k(E)=p_\ast(e(E)^{k+1})\in H^{2k}(B;\ZZ).
\]

The restriction of the Euler class $e(E)$ to
$\widetilde{B}\subseteq E$ gives another Euler class
$e(\widetilde{B})\in H^2(\widetilde{B};\ZZ)$. The push-down of the
$k$-th power of this second class is the \emph{boundary class}
\[
    \g_k(E)=\pi_\ast(e(\widetilde{B})^k)\in H^{2k}(B;\ZZ)
\]

The \emph{adjusted} or \emph{punctured} Miller-Morita-Mumford
classes $\ktilde_k(E)$ are given by
\[
    \ktilde_k(E)=\k_k(E)-\g_k(E).
\]

The surface bundle $\Sig_g^s\to E\rarrow B$ is classified by a map
$f:B\to BM_g^s$ and all three cohomology classes defined above are
pull-backs of universal classes
\[
    \k_k,\g_k,\ktilde_k\in H^{2k}(M_g^s;\ZZ)\cong H^{2k}(\Fat_g^s;\ZZ).
\]

By the fundamental results of Morita \cite{[Morita87]} and Miller
\cite{[Miller86:MMM]} the Miller-Morita-Mumford classes $\k_k$ and
the first $s$ boundary classes $\g_1,\cdots,\g_s$ are
algebraically independent over $\QQ$ in the stable range (given by
Harer stability \cite{[Harer85:Stability]}). The classes $\g_r$
for $r>s$ are polynomials in $\g_1,\cdots,\g_s$. For example,
\[
    \g_3=3\g_1\g_2-\frac12\g_1^3
\]
if $s=2$. Since the category $\Fat$ is the disjoint union of
$\Fat_g^s$ for all $g,s$, this also extends to $\ktilde_0=1-2g-s$
and we have the following.

\begin{thm}\label{thm:MMM classes are alg ind}
The universal adjusted Miller-Morita-Mumford classes
\[\ktilde_k\in H^\ast(\Fat;\QQ)\] for $k\geq0$ are algebraically
independent.
\end{thm}

In \cite{[I:BookOne]}, a combinatorial rational cocycle is
constructed for the punctured class $\ktilde_k$. It is called the
``cyclic set cocycle.''

\subsection{Cyclic set cocycle}

Let $\Z$ be the category of cyclically ordered sets and cyclic
order preserving monomorphisms. Then it is well known that
\[
    |\Z|\simeq BU(1).
\]
The $k$-th power of the first Chern class of the canonical complex
line bundle over $\Z$ is given by the \emph{unadjusted cyclic set
cocycle} $c_\Z^k$ whose value on a $2k$-simplex
\[
    C_\ast=(C_0\to C_1\to\cdots\to C_{2k})
\]
is given by
\begin{equation}\label{eq:cyclic set cocycle cZk}
    c_\Z^k(C_\ast)=\frac{\sum\sgn(a_0,a_1,\cdots,a_{2k})}
    {(-2)^k(2k-1)!!|C_0|\cdots|C_{2k}|}
\end{equation}
where the sum is taken over all choices of elements $a_i\in
C_i-C_{i-1}$. (In our notation we pretend that the maps in
$C_\ast$ are inclusion maps. Strictly speaking, $a_i$ should be an
element of $C_{2k}$ which lies in the image of $C_i$ but not in
the image of $C_{i-1}$.)

The unadjusted cyclic set cocycle has the following obvious
property.

\begin{prop}\label{prop:nec condition for cZk to be nonzero}
In order for the unadjusted cyclic set cocycle $c_\Z^k$ to be
nonzero on $C_\ast$ it is necessary (but not sufficient) for each
$C_i$ to be larger than $C_{i-1}$, i.e.,
\[
    |C_0|<|C_1|<\cdots<|C_{2k}|.
\]
\end{prop}

In \cite{[I:BookOne]} it is shown that the adjusted
Miller-Morita-Mumford class $\ktilde_k$ is given on the category
of ribbon graphs by evaluating the unadjusted cyclic set cocycle
(\ref{eq:cyclic set cocycle cZk}) at each vertex and dividing by
$-2$. Since there is already a factor of $(-2)^k$ in the
denominator it seems reasonable to define the (adjusted)
\emph{cyclic set cocycle} $\widetilde{c}_\Z^k$ by
\begin{equation}\label{eq:adjusted cyclic set cocycle on Z}
    \widetilde{c}_\Z^k(C_\ast)=\frac1{-2}c_\Z^k(C_\ast)=
    \frac{\sum\sgn(a_0,a_1,\cdots,a_{2k})}
    {(-2)^{k+1}(2k-1)!!|C_0|\cdots|C_{2k}|}.
\end{equation}

\begin{thm}[\cite{[I:BookOne]}]\label{thm:MMM is adj cyclic set cocycle}
The adjusted rational Miller-Morita-Mumford class $\ktilde_k\in
H^{2k}(\Fat;\QQ)$ is given on a $2k$-simplex
\begin{equation}\label{eq:2k simplex in Fat}
    \Gam_\ast=(\Gam_0\to\cdots\to \Gam_{2k})
\end{equation}
in $\Fat$ by evaluating the cyclic set cocycle
$\widetilde{c}_\Z^k$ on every vertex of $\Gam_0$ counted with
multiplicity.\footnote{The \emph{multiplicity} of a vertex is
defined to be the valence minus $2$.}
\end{thm}

We use the notation $\widetilde{c}_\Fat^{\,k}$ to denote this
rational cocycle on the cellular chain complex of $\Fat$ and we
refer to it by the same name, i.e., $\widetilde{c}_\Fat^{\,k}$ is
the \emph{cyclic set cocycle} on the category of ribbon graphs.

By definition, $\widetilde{c}_\Fat^{\,k}$ satisfies the following
important condition.

\begin{prop}\label{prop:cFatk is nonzero only when a vertex
increases in valence} The value of $\widetilde{c}_\Fat^{\,k}$ on a
(\ref{eq:2k simplex in Fat}) can be nonzero only if there is at
least one vertex $v_0$ of $\Gam_0$ whose image in each $\Gam_i$
has greater valence than its image in $\Gam_{i-1}$.
\end{prop}

Using dual cells we can pull back the adjusted
Miller-Morita-Mumford classes to the rational cohomology of the
graph cohomology complex $\G_\ast$.

\subsection{Adjusted MMM classes in $H^{2k}(\G_\ast;\QQ)$}

Recall that there is a chain homotopy equivalence
\[
    \psi:\G_\ast\otimes\QQ\to
    C_\ast(\Fat;\QQ)
\]
given up to sign by sending each rational generator $\<\Gam\>\in
\G_n$ to the average dual cell of $\Gam$ inside a finite model
$\thin$.
\[
    \psi\<\Gam\>=(-1)^{\binom{n+1}{2}}\overline{D}(\Gam).
\]
Since the cyclic set cocycle $\widetilde{c}_\Fat^{\,k}$ has the
same value on isomorphic simplices, its value on any dual cell of
$\Gam$ is the same. Therefore, its value of the average dual cell
$\overline{D}(\Gam)$ is equal to its value on any particular dual
cell $D(\Gam)$.

\begin{defn}
The \emph{cyclic set cocycle}
\[
    \ccc^k:\G_{2k}\to\QQ
\]
on the graph cohomology complex $\G_\ast$ is given by
\[
    \ccc^k\<\Gam\>=(-1)^k\widetilde{c}_\Fat^{\,k}(D(\Gam))
\]
for any choice of dual cell $D(\Gam)$.
\end{defn}

\begin{thm}
The cyclic set cocycle $\ccc^k$ represents the adjusted
Miller-Morita-Mumford class:
\[
    [\ccc^k]=\psi^\ast(\ktilde_k)\in H^{2k}(\G_\ast;\QQ).
\]
\end{thm}

It is clear from the definition of the cyclic set cocycle that it
can only be nonzero on the Witten cycle $W_{k}$ and it has the
same value on $\<\Gam\>$ for every element $\Gam$ of $W_{k}$ with
natural orientation. Therefore, $\ccc^k$ is proportional to the
dual Witten cycle $W_{k}^\ast$. The proportionality constant was
computed in \cite{[I:MMM_and_Witten]}.

\begin{thm}\label{thm:W2k in terms of adj ccc}
We have the following equation of rational cocycles on the
associative graph cohomology complex $\G_\ast$.
\[
    W_{k}^\ast=(-2)^{k+1}(2k+1)!!\ccc^k
\]
\end{thm}

\begin{proof}
As we showed in Theorem \ref{thm:orientation of simplices in
Keven} and Remark \ref{rem:orientation of K agrees with psi}, the
sign convention used in \cite{[I:MMM_and_Witten]} for every
simplex in the dual cell of any $\Gam\in W_{k}$ is given by
$(-1)^ko(\Gam_\ast)$ which agrees with the coefficient of the same
term in
\[
    \psi\<\Gam\>=(-1)^k\sum o(\Gam_\ast)(\Gam_0,\cdots,\Gam_{2k}).
\]
\end{proof}

Putting these together we get the following which can be
interpreted as statement about the relationship between rational
cohomology classes in graph cohomology, the category of ribbon
graphs or the mapping class group.

\begin{cor}[\cite{[I:MMM_and_Witten]}]\label{cor:Witten in terms of MMM}
The dual Witten cycle $[W_{n}^\ast]$ is a multiple of the adjusted
Miller-Morita-Mumford class:
\[
    [W_{n}^\ast]=(-2)^{n+1}(2n+1)!!\ktilde_n.
\]
\end{cor}

For example, we have:
\begin{align*}
    [W_0^\ast]&= -2\ktilde_0\\
    [W_1^\ast]&=\ 12\ktilde_1\\
    [W_2^\ast]&= -120\ktilde_2\\
    [W_3^\ast]&=\ 1680\ktilde_3\\
    [W_4^\ast]&= -30240\ktilde_4
\end{align*}

\subsection{Cup products of adjusted MMM classes}

We now consider cup products of the $\ktilde_k$'s for $k\geq1$. If
\[
    \ll=(\ll_1,\ll_2,\cdots,\ll_s)
\]
is a \emph{partition} of $n$ in the sense that $\ll_i>0$ and
$n=\sum \ll_i$ then let $\ktilde_{\ll}$ denote the cup product
\[
    \ktilde_{\ll_1}\cup\cdots\cup\ktilde_{\ll_s}\in H^{2n}(\Fat;\ZZ).
\]
Since the $\ktilde_{\ll_i}$ are even degree classes, this cup
product does not depend on the order of the $\ll_i$'s. However, at
the chain level, the order does matter. (The order matters in
$\Fat$ but not in the graph cohomology complex by Corollary
\ref{cor:ccc n and bnk are independent of order of n} below.) Let
$\widetilde{c}_\Fat^{\ll}$ denote the cup product
\begin{equation}\label{eq:cup product of cFats}
    \widetilde{c}_\Fat^{\ll}=
    \widetilde{c}_\Fat^{\ll_1}\cup\cdots\cup \widetilde{c}_\Fat^{\ll_s}\in
    C^{2n}(\Fat;\QQ)
\end{equation}
and let \[
    \ccc^{\ll}=\psi^\ast \widetilde{c}_\Fat^{\ll}\in \Hom(\G_{2n},\QQ)
\] denote the
pull-back of $\widetilde{c}_\Fat^{\ll}$ to $\G_\ast$. Then, by
Theorem \ref{thm:MMM is adj cyclic set cocycle}, we have
\[
    [\ccc^\ll]=\psi^\ast\ktilde_\ll.
\]

\begin{lem}\label{lem:cup products of ccc are only nonzero on W2n}
If $\ccc^{\ll}$ is nonzero on $\<\Gam\>\in \G_{2n}$ then $\Gam$
must lie in some Kontsevich cycle $W_{\mu}$ where $\ll$ is a
{refinement} of $\mu$.
\end{lem}

\begin{defn}\label{rem:pi represents nast as refinement of kast}
$\ll=(\ll_1,\cdots,\ll_s)$ is a \emph{refinement} of
$\mu=(\mu_1,\cdots,\mu_r)$ and we write
\[
    \ll\leq\mu
\]
if each $\mu_i$ is a sum of $\ll_j$'s so that, if we let $\pi(i)$
be the set of these indices $j$, then $\pi$ is a partition of the
set $\{1,2,\cdots,s\}$ into $r$ parts and
\[
    \ll_{\pi(i)}:=\sum_{j\in\pi(i)}\ll_j=\mu_i.
\]
We say that $\pi$ \emph{represents} $\ll$ as a refinement of
$\mu$.
\end{defn}

\begin{proof}
In order for the cup product (\ref{eq:cup product of cFats}) to be
nonzero on a dual cell $D(\Gam)$ where $\Gam\in\G_{2n}$ there must
be a \emph{nondegenerate} $2n$-simplex (i.e., where none of the
morphisms are isomorphisms)
\[
    \Gam_0\to\Gam_1\to\cdots\to\Gam_{2n}\cong\Gam
\]
so that
\[
    \widetilde{c}_\Fat^{\ll_1}(\Gam_0\to\cdots\to\Gam_{2\ll_1})\noteq0,
\]
\[
    \widetilde{c}_\Fat^{\ll_2}(\Gam_{2\ll_1}\to\cdots\to\Gam_{2\ll_1+2\ll_2})\noteq0,
    \text{ etc.}
\]
By Proposition \ref{prop:cFatk is nonzero only when a vertex
increases in valence}, $\Gam_{2\ll_1}$ must have a vertex of
multiplicity $2\ll_1$ and it must have a vertex which increases in
multiplicity by $2\ll_2$ by the time it gets to
$\Gam_{2\ll_1+2\ll_2}$. This implies that $\Gam_{2\ll_1+2\ll_2}$
must have either a vertex of multiplicity ${2\ll_1+2\ll_2}$ or two
vertices of multiplicity $2\ll_1, 2\ll_2$, resp.

By induction,
$\Gam=\Gam_{2\ll_1+\cdots+2\ll_{s-1}}=\Gam_{2n-2\ll_s}$ must lie
in a Kontsevich cycle $W_{\nu}$ so that $(\ll_1,\cdots,\ll_{s-1})$
is a refinement of $\nu$ as a partition of $n-\ll_s$. In order for
$\widetilde{c}_\Fat^{\,n-\ll_s}$ to be nonzero on the back $\ll_s$
face of $\Gam_\ast$, the graph $\Gam$ must have a vertex which
increases in multiplicity by $2\ll_s$ by the time it gets to
$\Gam_{2n}$. Thus $\Gam_{2n}$ must lie in $W_{\mu}$ where either
$\mu=(\nu,\ll_s)$ or $\mu$ is equal to $\nu$ with one of the
$\nu_i$ increased by $\ll_s$. In either cases,
\[
    \ll \leq (\nu,\ll_s)\leq \mu
\]
as claimed.
\end{proof}

\begin{lem}\label{lem:ccc n on W2k depends only on numbers}
Suppose that $\Gam\in W_{\mu}$ where $\mu=(\mu_1,\cdots,\mu_r)$ is
a partition of $n$ and $\ll$ is a {refinement} of $\mu$. Then the
value of $\ccc^{\ll}$ on $\<\Gam\>$ depends only on the ordered
partition $\ll$ and the unordered partition $\mu$ (and is
independent of the choice of $\Gam\in W_\mu$).
\end{lem}

\begin{rem}\label{rem:the ordering of the partitions n,k} We will
denote this number by
\begin{equation}\label{eq:the number b sub n upper k}
    b_{\ll}^{\mu}:=\ccc^{\ll}\<\Gam\>=\psi^\ast
    \widetilde{c}_\Fat^{\ll}\<\Gam\>
    =(-1)^n\widetilde{c}_\Fat^{\ll}D(\Gam)\in\QQ.
\end{equation}
It is obvious that the order of the $\mu_i$'s is not important. We
will later show (Corollary \ref{cor:ccc n and bnk are independent
of order of n}) that the order of the $\ll_j$'s is also not
important.
\end{rem}

\begin{proof}
If $\Gam$ lies in $W_{\mu}$ then $\Gam$ has vertices
$v_1,\cdots,v_r$ of codimension $2\mu_1,\cdots,2\mu_r$. In any
dual cell for $\Gam$, each of these vertices is expanded to a tree
in all possible ways (up to isomorphism). The rest of the graph is
left fixed. However, the cyclic set cocycle is only evaluated on
the vertices of these trees. Since the orientation of the
simplices in the dual cell depend only on these trees and the
value of the cocycle depends only on the trees, which in turn
depend only on the numbers $\mu_1,\cdots,\mu_r$ the value of
$\widetilde{c}_\Fat^{\ll}$ on $D(\Gam)$ depends only on $\mu$.
\end{proof}

Putting these two lemmas together we get the following.

\begin{thm}\label{thm:ccc is a sum of Wk's} Any cup product of
cyclic set cocycles
\[
    \ccc^{\ll}=\psi^\ast\left(\widetilde{c}_\Fat^{\ll_1}\cup
    \cdots\cup \widetilde{c}_\Fat^{\ll_s}\right)
\]
can be expressed as a rational linear combinations of dual
Kontsevich cycles by
\[
    \ccc^{\ll}=\sum b_{\ll}^{\mu}W_{\mu}^\ast
\]
where the sum is over all partitions $\mu$ of $n=\sum \ll_j$ so
that $\ll$ is a refinement of $\mu$ and $b_{\ll}^{\mu}$ is given
by (\ref{eq:the number b sub n upper k}) above.
\end{thm}

\begin{rem}\label{rem:we need to show Wmu are lin ind}
At the level of cohomology this theorem implies that
\[
    \ktilde_{\ll}=\sum b_{\ll}^{\mu}[W_{\mu}^\ast].
\]
Once we show that the cohomology classes $[W_\mu^\ast]$ are
linearly independent (Corollary \ref{cor:W star are lin ind}) then
this equation can be used to define $b_\ll^\mu$. We can then
conclude that $b_\ll^\mu$ and $\ccc^\ll$ are independent of the
order of $\ll$ (Corollary \ref{cor:ccc n and bnk are independent
of order of n}).
\end{rem}

\subsection{Computing the numbers $b_{\ll}^{\mu}$}

We first show that the computation of the numbers $b_{\ll}^{\mu}$
can be reduced to the case when $\mu=n$ is the trivial partition
of $n$.

\begin{lem}\label{lem:it suffices to work in An}
If $\mu=(\mu_1,\cdots,\mu_r)$ and $\ll=(\ll_1,\cdots,\ll_s)$ are
partitions of $n$ then
\begin{equation}\label{eq:b as sum of products}
    b_{\ll}^{\mu}=\sum_{\pi}\prod_{i=1}^r
    b_{\ll_{\pi(i)}}^{\mu_i}
\end{equation}
where we sum over all partitions $\pi$ of the set $\{1,\cdots,s\}$
into $r$ parts $\pi(i)$ representing $\ll$ as a refinement of
$\mu$ in the sense that $\mu_i=\sum_{j\in\pi(i)}\ll_j$.
\end{lem}

\begin{eg}\label{eg:b 32111 of 35}
Take the partition $\mu=(5,3)$ of $n=8$ and the refinement
$\ll=(3,2,1,1,1)$. Then
\[
    b_{3,2,1,1,1}^{5,3}=b_{3,2}^5b_{1,1,1}^3+3b_{3,1,1}^5b_{2,1}^3+b_{2,1,1,1}^5b_{3}^3.
\]
These terms come from the 5 ways in which $(5,3)$ can be refined
to $(3,2,1,1,1)$. They are
\[
    (\{3,2\},\{1,1,1\}), 3\times(\{3,1,1\},\{2,1\}), (\{2,1,1,1\},\{3\}).
\]
There are three ways to do the second refinement depending on
which $1$ goes to the right (to $3$ in $(5,3)$).
\end{eg}

\begin{proof} The number $b_{\ll}^{\mu}$ from (\ref{eq:the number b sub n upper
k}) times $(-1)^n$ is given by evaluating the cup product
$\widetilde{c}_\Fat^{\ll}$ of the cocycles
$\widetilde{c}_\Fat^{\ll_j}$ on every term of the dual cell
$D(\Gam)$ of any ribbon graph $\Gam$ in $W_{\mu}$.

The dual cell is given by choosing graphs over $\Gam$ and taking
nondegenerate $2n$-simplices
\[
    \Gam_\ast=(\Gam_0\to\cdots\to\Gam_{2n}\cong\Gam)
\]
(times $o(\Gam_\ast)$) where each $\Gam_i$ is a chosen
representative.

Look at all the terms in the dual cell $D(\Gam)$ which begin with
a fixed $\Gam_0$ (fixed as an object over $\Gam$). We will see
that the sum of the values of the cup product
$\widetilde{c}_\Fat^{\ll}$ on these terms is a sum of products
corresponding to the sum of products on the right hand side of
(\ref{eq:b as sum of products}).

Let $v_1,\cdots,v_r$ be the non-trivalent vertices of $\Gam$ and
let $T^1,\cdots,T^r$ be the trees in $\Gam_0$ which collapse to
these vertices. Thus $T^i$ has $2\mu_i$ edges which are naturally
ordered up to even permutation.

If the cup product
$\widetilde{c}_\Fat^{\ll}=\cup\widetilde{c}_\Fat^{\ll_p}$ is
nontrivial on a $2n$-simplex $\Gam_\ast$ beginning with $\Gam_0$
then to each index $p$ the cocycle $\widetilde{c}_\Fat^{\ll_p}$
must be nontrivial on
\begin{equation}\label{eq:middle np simplex}
    \Gam_{2n_p}\to\Gam_{2n_p+1}\to\cdots\to\Gam_{2n_p+2\ll_p}
\end{equation}
where $n_p=\ll_1+\ll_2+\cdots+\ll_{p-1}$. This means the edges
which collapse in this sequence must all lie in the same tree
$T^i$. Let $i=f(p)$. Then $f$ is an epimorphism
\begin{equation}\label{eq:the epimorphism f}
    f:\{1,\cdots,s\}\twoheadrightarrow \{1,\cdots,r\}
\end{equation}
so that $\mu_i$ is equal to the sum of the $\ll_p$ for all $p\in
f^{-1}(i)$. In other words, $\pi=f^{-1}$ represents $\ll$ as a
refinement of $\mu$.

The value of the cocycle $\widetilde{c}_\Fat^{\ll_p}$ on the
middle $2\ll_p$-simplex (\ref{eq:middle np simplex}) depends only
on the sequence of edges in $T^i$ which are collapsing. Thus, if
we fix the epimorphism (\ref{eq:the epimorphism f}), then, for
each index $i$, the order in which the edges of $T^i$ collapse in
(\ref{eq:middle np simplex}) varies independently of the edges of
$T^j$ for $j\noteq i$. Consequently, the sum of products becomes a
product of sums ($\pi=f^{-1}$ still being fixed).
\[
    \sum_{\Gam_\ast}\prod_{i=1}^r\prod_{p\in\pi(i)}
    \widetilde{c}_\Fat^{\ll_p}(\Gam_{2n_p},\cdots,\Gam_{2n_p+2\ll_p})
    =\prod_{i=1}^r\sum_{T_\ast^i}\prod_{p\in\pi(i)}
    \widetilde{c}_\Fat^{\ll_p}(T^i_{2n_p},\cdots,T^i_{2n_p+2\ll_p})
\]
where $T^i_k$ is the inverse image in $\Gam_k$ of the vertex $v_i$
of $\Gam$.

But a different example with $r=1$ gives the same sum:
\[
    (-1)^{\mu_i}b_{\ll_{pi(i)}}^{\mu_i}=\sum_{T_\ast^i}\prod_{p\in\pi(i)}
    \widetilde{c}_\Fat^{\ll_p}(T^i_{2n_p},\cdots,T^i_{2n_p+2\ll_p}).
\]
Taking the product over all $i$ and the sum over all $\pi$ we get
\[
        (-1)^n b_{\ll}^{\mu}=\sum_{\pi}\prod_{i=1}^r
    (-1)^{\mu_i}b_{\ll_{\pi(i)}}^{\mu_i}
\]
which is the same as (\ref{eq:b as sum of products}) since $n=\sum
\mu_i$.
\end{proof}

We can now compute some of the numbers $b_{\ll}^{\mu}$. We start
with the following case which follows from Theorem \ref{thm:W2k in
terms of adj ccc}.

\begin{lem}\label{lem:the value of b kk} In the case of the
trivial partitions $\ll=\mu=n$ we have:
\[
    b_n^n=\frac1{(-2)^{n+1}(2n+1)!!}.
\]
\end{lem}

By Lemma \ref{lem:it suffices to work in An} this give the
following.

\begin{prop}\label{prop:value of b kstar kstar}
If the partition $\ll$ of $n$ has $m_i$ terms equal to $\ll_i$ for
$i=1,\cdots,r$ (so that $\sum m_i\ll_i=n$) then
\[
    b_{\ll}^{\ll}=
    \prod_{i=1}^r{m_i!}\left(
    b_{\ll_i}^{\ll_i}
    \right)^{m_i}
    =
    \prod_{i=1}^r\frac{m_i!}
    {\left((-2)^{\ll_i+1}(2\ll_i+1)!!\right)^{m_i}}.
\]
\end{prop}

\begin{proof}
The factor of $\prod m_i!$ is equal to the number of ways that the
partition $\ll$ refines itself.
\end{proof}

\subsection{Kontsevich cycles in terms of MMM classes}

We are now ready to show that the dual Konsevich cycles
$W_{\mu}^\ast$ represent polynomials in the adjusted
Miller-Morita-Mumford classes $\ktilde_k$ (as we claimed in
\cite{[I:MMM_and_Witten]}). We will then conclude that their
cohomology classes $[W_\mu^\ast]$ are linearly independent as
promised in Remark \ref{rem:we need to show Wmu are lin ind}.

To avoid circular reasoning, we must assume at this point that the
numbers $b_\ll^\mu$ may depend on the order of the parts of $\ll$.
Therefore, we take the ordering of both $\ll$ and $\mu$ to be
nonincreasing: $\ll_1\geq\ll_2\geq\cdots\geq\ll_r>0$. If $\ll$ is
a refinement of $\mu$ then $\ll\leq\mu$ in lexicographic order.
Consequently, the matrix
\[
    B_n=(b_\ll^\mu)
\]
is upper triangular. (A priori this uses only some of the numbers
$b_\ll^\mu$.) The diagonal entries $b_\ll^\ll$ are nonzero by
Proposition \ref{prop:value of b kstar kstar} so $B_n$ is
invertible. Let $A_n$ be the inverse matrix
\[
    A_n=B_n^{-1}=(a_\mu^\ll).
\]
The entries of this matrix are rational numbers uniquely
determined by the equation
\begin{equation}\label{eq:defining eq for a ll mu}
    \sum_{\nu}a_\ll^\nu b_\nu^\mu=\delta_\ll^\mu.
\end{equation}

The main theorem is the following.

\begin{thm}\label{thm:W star is a polynomial in MMM}
The cohomology classes of the dual Kontsevich cycles are
polynomials in the adjusted Miller-Morita-Mumford classes:
\[
    [W_{\mu}^\ast]=\sum_{\ll}a_{\mu}^{\ll}\ktilde_{\ll}.
\]
\end{thm}

\begin{rem}\label{rem:W star in terms of MMM holds unstably}
This formula holds in the rational cohomology of $\G_\ast$, and in
the integral cohomology modulo torsion of the category of ribbon
graphs and the mapping class group $M_g^s$ for all $g,s\geq1$ (and
for $g=0,s\geq3$).
\end{rem}

\begin{proof}
Since the matrix $B_n$ is invertible, the system of linear
equations given in Theorem \ref{thm:ccc is a sum of Wk's} has the
unique solution
\[
    W_{\mu}^\ast=
    \sum_{\ll}a_{\mu}^{\ll}\ccc^{\ll}
\]
provided that $\ll$ is in nonincreasing order.
\end{proof}

Since the coefficients $a_{\mu}^{\ll}$ form an invertible matrix
and the monomials $\ktilde_{\ll}$ are linearly independent by
Theorem \ref{thm:MMM classes are alg ind} we have the following.

\begin{cor}\label{cor:W star are lin ind}
The cohomology classes $[W_{\mu}^\ast]$ are linearly independent
over $\QQ$.
\end{cor}

This, in turn, implies the following as explained in Remark
\ref{rem:we need to show Wmu are lin ind}.

\begin{cor}\label{cor:ccc n and bnk are independent of order of n}
The cocycles $\ccc^{\ll}$ and the numbers $b_{\ll}^{\mu}$ are
independent of the ordering of the partitions $\ll$, $\mu$.
\end{cor}

Combining these with Theorem \ref{thm:ccc is a sum of Wk's}, we
get the following.

\begin{cor}\label{cor:span W ll = QQ[k1,k2,etc]}
The cohomology classes $[W_\mu^\ast]$ form a $\QQ$-basis for the
polynomial algebra generated by the adjust Miller-Morita-Mumford
classes $\ktilde_k$ for $k\geq1$.
\end{cor}

\subsection{Computing $a_\ll^\mu$}

We will use the defining equation (\ref{eq:defining eq for a ll
mu}) to determine the numbers $a_\ll^\mu$ in some simple cases.

\begin{prop}\label{prop:a nn}
When $\ll=\mu=n$ we have
\[
    a_n^n=\frac1{b_n^n}=(-2)^{n+1}(2n+1)!!
\]
\end{prop}

To simplify the notation we write $a_n$ and $b_n$ instead of
$a_n^n$ and $b_n^n$.

\begin{prop}\label{prop:value of a ll ll}
When $\ll=\mu=(\ll_1^{m_1},\cdots,\ll_r^{m_r})$ with $\sum
m_i\ll_i=n$ then
\[
    a_\ll^\ll=\frac1{b_\ll^\ll}=\prod_{i=1}^r
    \frac{\left((-2)^{k_i+1}(2k_i+1)!!\right)^{m_i}}{m_i!}.
\]
\end{prop}

For example, $a_{2,2}^{2,2}=120^2/2=7200$ and
$a_{1,1,1}^{1,1,1}=12^3/3!=288$. These numbers are not always
integers. For example,
\[
    a_{1^5}^{1^5}=\frac{12^5}{5!}=\frac{10368}{5}.
\]

\begin{cor}\label{cor:W ast is a polynomial in MMM}
The dual Kontsevich cycle $[W_{\ll}^\ast]$ is a polynomial of
degree $\sum m_i$ in the adjusted Miller-Morita-Mumford classes
with leading term
\[
    \prod_{i=1}^r\frac{\left((-2)^{\ll_i+1}(2\ll_i+1)!!\ktilde_{\ll_i}\right)^{m_i}}
{m_i!}
\]
\end{cor}

In order to compute the remaining terms we need to compute
$a_\ll^\mu$ for $\ll\noteq\mu$. The first case is $a_{n,k}^{n+k}$
which occurs in the equation
\begin{equation}\label{eq:W 2n2k}
    [W_{n,k}^\ast]=a_{n,k}^{n,k}\ktilde_n\ktilde_k+a_{n,k}^{n+k}\ktilde_{n+k}.
\end{equation}
 The defining equation for $a_{n,k}^{n+k}$ is
\[
    a_{n,k}^{n+k} b_{n+k}+a_{n,k}^{n,k}b_{n,k}^{n+k}=0.
\]
From this and Proposition \ref{prop:a nn} we get:
\begin{align*}
    a_{n,k}^{n,k}&=(-2)^{2n+2k+2}(2n+1)!!(2k+1)!!\\
    a_{n,k}^{n+k}&=-a_{n+k}a_{n,k}^{n,k}b_{n,k}^{n+k}\\
    &=-(-2)^{2n+2k+3}(2n+2k+1)!!(2n+1)!!(2k+1)!!b_{n,k}^{n+k}
\end{align*}
if $n\noteq k$ and
\begin{align*}
    a_{n,n}^{n,n}&=-(-2)^{4n+3}\left((2n+1)!!\right)^2\\
    a_{n,n}^{2n}&=(-2)^{4n+2}(4n+1)!!\left((2n+1)!!\right)^2b_{n,n}^{2n}
\end{align*}

In the next section we will compute $b_{n,k}^{n+k}$ in the special
case $k=1$.

% \vfill\eject

\setcounter{section}{3}

\section{Some computations} % sec 04

\begin{enumerate}
    \item {The degenerate case $n=0$}
    \item {Computation of $b_{n,1}^{n+1}$}
    \item {Conjectures}
\end{enumerate}

We will compute the numbers $b_{n,1}^{n+1}$ and obtain
$a_{n,1}^{n+1}$ and the coefficients of $[W_{n,1}^\ast]$ when
expressed as a polynomial in the adjusted Miller-Morita-Mumford
classes.

\subsection{The degenerate case $n=0$}

First we consider the degenerate case $n=0$. In this case it is
easy to compute $W_{k,0}^\ast$ and work backwards since this
cocycle counts the number of pairs of vertices, one of codimension
$2k$ and the other of codimension 0 (i.e., trivalent). The number
of trivalent vertices of a graph with one vertex of multiplicity
$2k+1$ is $t$ where
\[
    t+2k+1=-2\chi=-2\ktilde_0.
\]
So, $t=-2\ktilde_0-(2k+1)$ and
\[
    [W_{k,0}^\ast]=t[W_{k}^\ast]=(-2\ktilde_0-(2k+1))(-2)^{k+1}(2k+1)!!\ktilde_k
\]
\begin{equation}\label{eq:Wk0}
    [W_{k,0}^\ast]=(-2)^{k+2}(2k+1)!!\ktilde_k\ktilde_0-(2k+1)(-2)^{k+1}(2k+1)!!\ktilde_k.
    \end{equation}
    The right hand side should be divided by $2$ when $k=0$.

For the case $k=1$ this gives
\[
    [W_{1,0}^\ast]=[W_{0,1}^\ast]=-24\ktilde_0\ktilde_1-36\ktilde_1.
\]
Consequently, $a_{0,1}^1=-36$ and
\[
    b_{0,1}^1=\frac1{2^59}a_{0,1}^1=-\frac18.
\]

\subsection{Computation of $b_{n,1}^{n+1}$}

To obtain $b_{n,1}^{n+1}$ in general we use the formula from Lemma
\ref{lem:ccc n on W2k depends only on numbers}:
\[
    b_{n,1}^{n+1}=(-1)^{n+1}(\widetilde{c}_\Fat^{\,n}\cup
    \widetilde{c}_\Fat^{\,1})D(\Gam)
\]
where $\Gam$ is any graph in the Witten cycle $W_{n+1}$, i.e., a
ribbon graph which is trivalent except at one vertex of
multiplicity $2n+2$.

The dual cell $D(\Gam)$ is a sum
\[
    D(\Gam)=\sum_{\Gam_\ast}o(\Gam_\ast)(\Gam_0,\cdots,\Gam_{2n+2}=\Gam).
\]
When we evaluate the cup product
$(-1)^{n+1}(\widetilde{c}_\Fat^{\,n}\cup
    \widetilde{c}_\Fat^{\,1})$ we get
    \[
b_{n,1}^{n+1}=(-1)^{n+1}\sum_{\Gam_\ast}o(\Gam_\ast)
\widetilde{c}_\Fat^{\,n}(\Gam_0,\cdots,\Gam_{2n})
\widetilde{c}_\Fat^{\,1}(\Gam_{2n},\Gam_{2n+1},\Gam)
    \]

There are three possible configurations for $\Gam_{2n}$. They are
Case 1, given in Figure \ref{fig:tree01} and Case 2a,2b, given in
Figure \ref{fig:tree02} from section 2. In all cases, $\Gam_{2n}$
is odd-valent so it has a natural orientation. The orientation of
the sequence $\Gam_\ast$ can then be expressed as a product of two
orientations:
\[
    o(\Gam_0,\cdots,\Gam_{2n},\Gam',\Gam)=o(\Gam_0,\cdots,\Gam_{2n})
    o(\Gam_{2n},\Gam',\Gam)
\]
where we write $\Gam'=\Gam_{2n+1}$. Consequently, for any fixed
$\Gam_{2n}$, we get
\[
(-1)^n\sum o(\Gam_0,\cdots,\Gam_{2n})\widetilde{c}_\Fat^{\,n}
(\Gam_0,\cdots,\Gam_{2n})
    \left(
-\sum_{\Gam'}
o(\Gam_{2n},\Gam',\Gam)
\widetilde{c}_\Fat^{\,1}(\Gam_{2n},\Gam',\Gam)
    \right).
\]
The first factor is $b_n$ regardless of $\Gam_{2n}$ so
\[
    \frac{b_{n,1}^{n+1}}{b_n}=-\sum_{\Gam_{2n},\Gam'}
    o(\Gam_{2n},\Gam',\Gam)\widetilde{c}_\Fat^{\,1}(\Gam_{2n},\Gam',\Gam).
\]
If we denote the sequence $\Gam_{2n}\to\Gam'\to\Gam$ by
$\Gam_\ast$ and substitute $1/b_n=a_n$ we get
\begin{equation}\label{eq:sum for an bn1 n+1}
    a_nb_{n,1}^{n+1}=-\sum_{\Gam_\ast}o(\Gam_\ast)
    \widetilde{c}_\Fat^{\,1}(\Gam_\ast).
\end{equation}

Recall that
\[
    \widetilde{c}_\Fat^{\,1}(\Gam_\ast)=\sum
    \frac{\mu(v_0)}{4}\frac{\sum \sgn(a,b,c)}{|C_0||C_1||C_2|}
\]
with the first sum being over all vertices $v_0$ of $\Gam_{2n}$
where $\mu(v_0)$ is the multiplicity of $v_0$ and $|C_i|$ is the
valence of the image of $v_i$ in $\Gam_{2n+i}$. The sign sum
\[
    \sum\sgn(a,b,c)=\sgn(x_1\cdots x_{2n+5})
\]
is the number of times that the letters $a,b,c$ occur in the
correct cyclic order in the word $w=x_1\cdots x_{2n+5}$ minus the
number of times it occurs in the other cyclic order in $w$ where
the $j$th letter $x_j$ of $w$ is equal to $b_i$
($b_0=a,b_1=b,b_2=c$) if the $j$th region in the complement of the
graph reaches $v_0$ at the $i$th step ($\Gam_{2n+i}$).

The sum (\ref{eq:sum for an bn1 n+1}) breaks up into three parts
depending on the graph $\Gam_{2n}$.

\noi\underline{Case 1}. Suppose the graph $\Gam_{2n}$ is given by
Figure \ref{fig:tree01} and $T'=T_{2n}/e_1$. Then
\[
    o(\Gam_\ast)=(-1)^m
\]
The number of times this same configuration (with fixed $1\leq
m\leq2n+2$) occurs is
\[
    2n+5.
\]
The cyclic set cocycle $\widetilde{c}_\Fat^{\,1}(\Gam_\ast)$ has
two terms
\begin{enumerate}
    \item The center vertex has multiplicity $2n+1$ and average
    sign
    \[
\frac{\sgn(a^mca^{2n-m+3}b)}{(2n+3)(2n+4)(2n+5)} =
\frac{2n-2m+3}{(2n+3)(2n+4)(2n+5)}
    \] for a contribution of
    \[
\frac{(2n+1)(2n-2m+3)}{4(2n+3)(2n+4)(2n+5)}
    \]
    \item The vertex $v_1$ has multiplicity 1 and average sign
    \[
\frac{\sgn(a^3b^{m-1}cb^{2n-m+2})}{3(2n+4)(2n+5)} =
-\frac{2n-2m+3}{(2n+4)(2n+5)}
    \]
\end{enumerate}
The total value of the cocycle is
\[
    \widetilde{c}_\Fat^{\,1}(\Gam_\ast)=
    \frac{(2n-2m+3)(2n+1-(2n+3))}
    {4(2n+3)(2n+4)(2n+5)}
    =\frac{-(2n-2m+3)}{2(2n+3)(2n+4)(2n+5)}.
\]

Multiply this by $(-1)^m(2n+5)$ and sum over all $1\leq m\leq2n+2$
to get
\[
    \frac{-1}{2(2n+3)(2n+4)}\sum_{m=1}^{2n+2}(-1)^m
    (2n-2m+3)=\frac{n+1}{(2n+3)(2n+4)}.
\]

\noi\underline{Case 2a}. Suppose $\Gam_{2n}$ is given by Figure
\ref{fig:tree02}a and $\Gam'=\Gam_{2n}/e_1$. Then
\[
    o(\Gam_\ast)=\sgn(a_1,\cdots,a_{2n+3},b_1,b_2)=1.
\]
This configuration occurs $2n+5$ times and the cyclic set cocycle
has two terms.
\begin{enumerate}
    \item The center vertex has multiplicity $2n+1$ and average
    sign
    \[
\frac{\sgn(a^{2n+3}bc)}{(2n+3)(2n+4)(2n+5)}=\frac1{(2n+4)(2n+5)}
    \]
    for a contribution of
    \[
\frac{2n+1}{4(2n+4)(2n+5)}.
    \]
    \item The vertex $v_1$ has multiplicity 1 and average sign
    \[
\frac{\sgn(a^2cab^{2n+1})}{3(2n+4)(2n+5)}=\frac{-(2n+1)}{3(2n+4)(2n+5)}
    \]
\end{enumerate}
So
\[
    \widetilde{c}_\Fat^{\,1}(\Gam_\ast)=
    \frac{2n+1}{6(2n+4)(2n+5)}.
\]
Multiply this by $2n+5$ for a total of
\[
    \frac{2n+1}{6(2n+4)}.
\]

\noi\underline{Case 2a'}. Suppose that $T_{2n}$ is the same
(Figure \ref{fig:tree02}a) but $\Gam'=\Gam_{2n}/e_2$. Then
$o(\Gam_\ast)=-1$. The configuration still occurs $2n+5$ times and
the cyclic set cocycle again has two terms.
\begin{enumerate}
    \item At $v_1$ we have
    \[
    \frac14
\frac{\sgn(a^2bac^{2n+1})}{3\cdot4(2n+5)}= \frac{2n+1}{48(2n+5)}.
    \]
    \item At $v_2$ we get
    \[
    \frac14
\frac{\sgn(ba^3c^{2n+1})}{3\cdot4(2n+5)}=
\frac{-3(2n+1)}{48(2n+5)}.
    \]
\end{enumerate}
So,
\[
    \widetilde{c}_\Fat^{\,1}(\Gam_\ast)
    =\frac{-(2n+1)}{24(2n+5)}
\]
making the total in this case
\[
    \frac{2n+1}{24}.
\]

Case 2b is the same as Case 2a. (The sign changes twice.) So
\[
    a_nb_{n,1}^{n+1}=\frac{n+1}{(2n+3)(2n+4)}
    +\frac{2n+1}{3(2n+4)}+\frac{2n+1}{12}.
\]
Simplifying this expression we get:
\[
    a_nb_{n,1}^{n+1}=\frac{2n+5}{12}-\frac1{2(2n+3)}.
\]
Multiplying by $a_1a_{n+1}=12(-2)^{n+2}(2n+3)!!$ we get
\[
    a_{n+1}a_1a_nb_{n,1}^{n+1}=(-2)^{n+2}(2n+5)!!+3(-2)^{n+3}(2n+1)!!
\]

\begin{thm}\label{thm:W 2n2 equals}
For $n\noteq1$ we have
\[
    [W_{n,1}^\ast]=3(-2)^{n+3}(2n+1)!!(\ktilde_n\ktilde_1-\ktilde_{n+1})-
    (-2)^{n+2}(2n+5)!!\ktilde_{n+1}.
\]
For $n=1$ we divide the right hand side by $2$.
\end{thm}

For example, we have
\begin{align*}
    [W_{0,1}^\ast]&=-24\ktilde_0\ktilde_1-36\ktilde_1\\
    [W_{1,1}^\ast]&=\ 72\ktilde_1^2+348\ktilde_2\\
    [W_{2,1}^\ast]&=-1440\ktilde_2\ktilde_1-13680\ktilde_3\\
    [W_{3,1}^\ast]&=\ 20160\ktilde_3\ktilde_1+312480\ktilde_4\\
    [W_{4,1}^\ast]&=-362880\ktilde_4\ktilde_1-8285760\ktilde_5
\end{align*}

This agrees with the calculation (\ref{eq:Wk0}) when $n=0,k=1$ and
also agrees with the calculation of Arbarello and Cornalba
\cite{[Arbarello-Cornalba:96]}. The sign difference comes from the
fact that they use the opposite sign for all $\ktilde_{even}$.
(What we call $\ktilde_k$ is what they would call $(-1)^{k+1}\k_k$
restricted to the open moduli space of curves.) The calculations
of \cite{[Arbarello-Cornalba:96]} give cohomology classes which
they showed act as Poincar\'{e} duals of the Kontsevich cycles
with respect to products of boundary cycles. Therefore, Theorem
\ref{thm:MMM classes are alg ind} and Theorem \ref{thm:W star is a
polynomial in MMM} already suggests that they must be correct.

\subsection{Conjectures}

The formula in Theorem \ref{thm:W 2n2 equals} has an apparent
symmetry which also appears in the formula (\ref{eq:Wk0}) when it
is rephrased as follows.
\begin{equation}\label{eq:Wk0 rephrased}
    [W_{k,0}^\ast]=(-2)^{k+2}(2k+1)!!(\ktilde_k\ktilde_0-\ktilde_k)
    -(-2)^{k+1}(2k+3)!!\ktilde_k.
    \end{equation}
This leads to the following conjecture.

\begin{conj}\label{conj:W 2k 2n should be}
\[
[W_{n,k}^\ast]=(-2)^{n+k+2}(2n+1)!!(2k+1)!!(\ktilde_n\ktilde_k-\ktilde_{n+k})
    -(-2)^{n+k+1}(2n+2k+3)!!\ktilde_{n+k}.
\]
The right hand side should be divided by $2$ if $n=k$.
\end{conj}

\begin{rem}\label{rem:simplification of W 2n2k}
First, note that this conjectured formula can be simplified using
the numbers $a_n=(-2)^{n+1}(2n+1)!!$:
\[
    [W_{n,k}^\ast]=a_na_k(\ktilde_n\ktilde_k-\ktilde_{n+k})
    +\frac12 a_{n+k+1}\ktilde_{n+k}.
\]
Next, we also note that it is symmetrical in $n,k$. And finally,
in the first new case when $n=k=2$ it gives
\begin{equation}\label{eq:W22 calculation}
    [W_{2,2}^\ast]
    =\frac{120^2}{2}(\ktilde_2^2-\ktilde_4)
    +\frac{665280}{4}\ktilde_4
    =7200\ktilde_2^2+159120\ktilde_4
\end{equation}
which agrees with Arbarello and Cornalba
\cite{[Arbarello-Cornalba:96]}.
\end{rem}

\begin{rem}\label{rem:joint work with Kleber}
Michael Kleber and I have made further progress on calculation and
interpretation of the coefficients $b_\ll^n$. So far we verified
Conjecture \ref{conj:W 2k 2n should be} for all $k \leq7$. In
particular this proves (\ref{eq:W22 calculation}). We also found
that
\[
    [W_{1,1,1}^\ast]=288\ktilde_1^3+4176\ktilde_1\ktilde_2+20736\ktilde_3.
\]
This, together with Theorem \ref{thm:W 2n2 equals} and Corollary
\ref{cor:Witten in terms of MMM}, verifies all calculations of the
coefficients $a_\ll^\mu$ given by Arbarello and Cornalba. Details
will be given in a subsequent joint paper.
\end{rem}

%\vfill\eject

%%%%%%%%%%%%%%%%%%%%%%%%%%%%%%%%%%%%%%%%%%%%%%%%%%%%%%%%%%%%%%%%%%%%%%%%%%%%%%%%%%%%%%%%

\bibliographystyle{amsalpha}
\bibliography{C:/bookbib3}

\end{document}